\documentclass[onecolumn,11pt]{article}
\usepackage[top=1in, bottom=1in, left=1in, right=1in]{geometry}
\usepackage{amsmath,amsfonts,amscd,amssymb}
\usepackage{graphicx}
\usepackage{epstopdf}
\usepackage{overpic}
\usepackage{cancel}
\usepackage{rotating}
\usepackage{url}
\usepackage{caption}
\usepackage{color}
\usepackage{rotating}
\usepackage{multirow}
\usepackage{wrapfig}
\usepackage{mathtools}
\usepackage{subeqnarray}
\usepackage{setspace}
\usepackage{palatino} 
\usepackage[bottom,flushmargin,hang,multiple]{footmisc}
\usepackage{lipsum}
\newcommand\blfootnote[1]{%
  \begingroup
  \renewcommand\thefootnote{}\footnote{#1}%
  \addtocounter{footnote}{-1}%
  \endgroup
}

\newcommand{\bx}{\mathbf{x}}

\newcommand{\bX}{\mathbf{X}}

\newcommand{\bXi}{\boldsymbol{\Xi}}
\newcommand{\bxi}{\boldsymbol{\xi}}

\newcommand{\bTheta}{\boldsymbol{\Theta}}

\newcommand{\koop}{\mathcal{K}}

\title{\LARGE{\vspace{-.5in}\textbf{Chaos as an Intermittently Forced Linear System}}\vspace{-.15in}}
\author{\normalsize{Steven L. Brunton$^{1*}$, Bingni W. Brunton$^2$, Joshua L. Proctor$^3$, Eurika Kaiser$^{1}$, J. Nathan Kutz$^4$}\\
\footnotesize{$^1$ Department of Mechanical Engineering, University of Washington, Seattle, WA 98195, United States}\\
\footnotesize{$^2$ Department of Biology, University of Washington, Seattle, WA 98195, United States}\\
\footnotesize{$^3$Institute for Disease Modeling, Bellevue, WA 98004, United States}\\ 
\footnotesize{$^4$ Department of Applied Mathematics, University of Washington, Seattle, WA 98195, United States\vspace{-.1in}}
}
\date{}
\begin{document}
\maketitle

\blfootnote{$^*$ Corresponding author (sbrunton@uw.edu).\\ \noindent \textbf{Matlab code:}  http://faculty.washington.edu/sbrunton/HAVOK.zip\\
\noindent \textbf{Video abstract:}  http://youtu.be/831Ell3QNck}
\vspace{-.2in}
\begin{abstract}
Understanding the interplay of order and disorder in chaotic systems is a central challenge in modern quantitative science. 
We present a universal, data-driven decomposition of chaos as an intermittently forced linear system.  
This work combines Takens' delay embedding with modern Koopman operator theory and sparse regression to obtain linear representations of strongly nonlinear dynamics.  
The result is a decomposition of chaotic dynamics into a linear model in the leading delay coordinates with forcing by low energy delay coordinates; we call this the Hankel alternative view of Koopman (HAVOK) analysis.    
This analysis is applied to the canonical Lorenz system, 
as well as to real-world examples such as the Earth's magnetic field reversal, and data from electrocardiogram, electroencephalogram, and measles outbreaks.  
In each case, the forcing statistics are non-Gaussian, with long tails corresponding to rare events that trigger intermittent switching and bursting phenomena; this forcing is highly predictive, providing a clear signature that precedes these events.    
Moreover, the activity of the forcing signal demarcates large coherent regions of phase space where the dynamics are approximately linear from those that are strongly nonlinear.  \\

\end{abstract}


\vspace{-.15in}
\section{Introduction}
Dynamical systems describe the world around us, modeling the interactions between quantities that co-evolve in time~\cite{guckenheimer_holmes}.  
These dynamics often give rise to rich and complex behaviors that may be difficult to predict from uncertain measurements, a phenomena that is commonly known as \emph{chaos}.  
Chaotic dynamics are ubiquitous in the physical, biological, and engineering sciences, and they have captivated amateurs and experts alike for over a century.  
The motion of planets~\cite{Poincare1890am}, weather and climate~\cite{Lorenz1963jas,Majda2007bpnas,Majda2012nonlinearity,Giannakis2012pnas,Sapsis2013pnas,Majda2014pnas}, population dynamics~\cite{Bjornstad2001science,Sugihara2012science,Ye2015pnas}, epidemiology~\cite{Sugihara1990nature}, financial markets, earthquakes, solar flares, and turbulence~\cite{Kolmogorov1941b,Kolmogorov1941c,Kolmogorov1941,Takens1981lnm,Brunton2015amr}, all provide compelling examples of chaos.  
Despite the name, chaos is not random, but is instead highly organized, exhibiting coherent structure and patterns~\cite{Tsonis1992nature,Crutchfield2012naturephys}.  

The confluence of big data and advanced algorithms in machine learning is driving a paradigm shift in the analysis and understanding of dynamical systems in science and engineering.  
Data are abundant, while physical laws or governing equations remain elusive, as is true for problems in climate science, finance, and neuroscience.  
Even in classical fields such as turbulence, where governing equations do exist, researchers are increasingly turning towards data-driven analysis~\cite{Sapsis2013pnas,Majda2014pnas,Brunton2015amr}.  
Many critical data-driven problems, such as predicting climate change, understanding cognition from neural recordings, or controlling turbulence for energy efficient power production and transportation, are primed to take advantage of progress in the data-driven discovery of dynamics~\cite{Schmidt2009science,Brunton2016pnas}.

An early success of data-driven dynamical systems is the celebrated Takens embedding theorem~\cite{Takens1981lnm}, which allows for the reconstruction of an attractor that is diffeomorphic to the original chaotic attractor from a time series of a single measurement.  
This remarkable result states that, under certain conditions, the full dynamics of a system as complicated as a turbulent fluid may be uncovered from a time series of a single point measurement.  
Delay embedding has been widely used to analyze and characterize chaotic systems~\cite{Farmer1987prl,Crutchfield1987cs,Sugihara1990nature,Rowlands1992physD,Abarbanel1993rmp,Sugihara2012science,Ye2015pnas}, as well as for linear system identification with the eigensystem realization algorithm (ERA)~\cite{ERA:1985} and in climate science with the singular spectrum analysis (SSA)~\cite{Broomhead1989prsla} and nonlinear Laplacian spectrum analysis (NLSA)~\cite{Giannakis2012pnas}.  
Until now, there has been a disconnect between the use of delay embeddings to characterize chaos and their rigorous use to identify models of the nonlinear dynamics. 

\begin{figure}[b!]
\begin{center}
\vspace{-.1in}
\begin{overpic}[width=.9\textwidth]{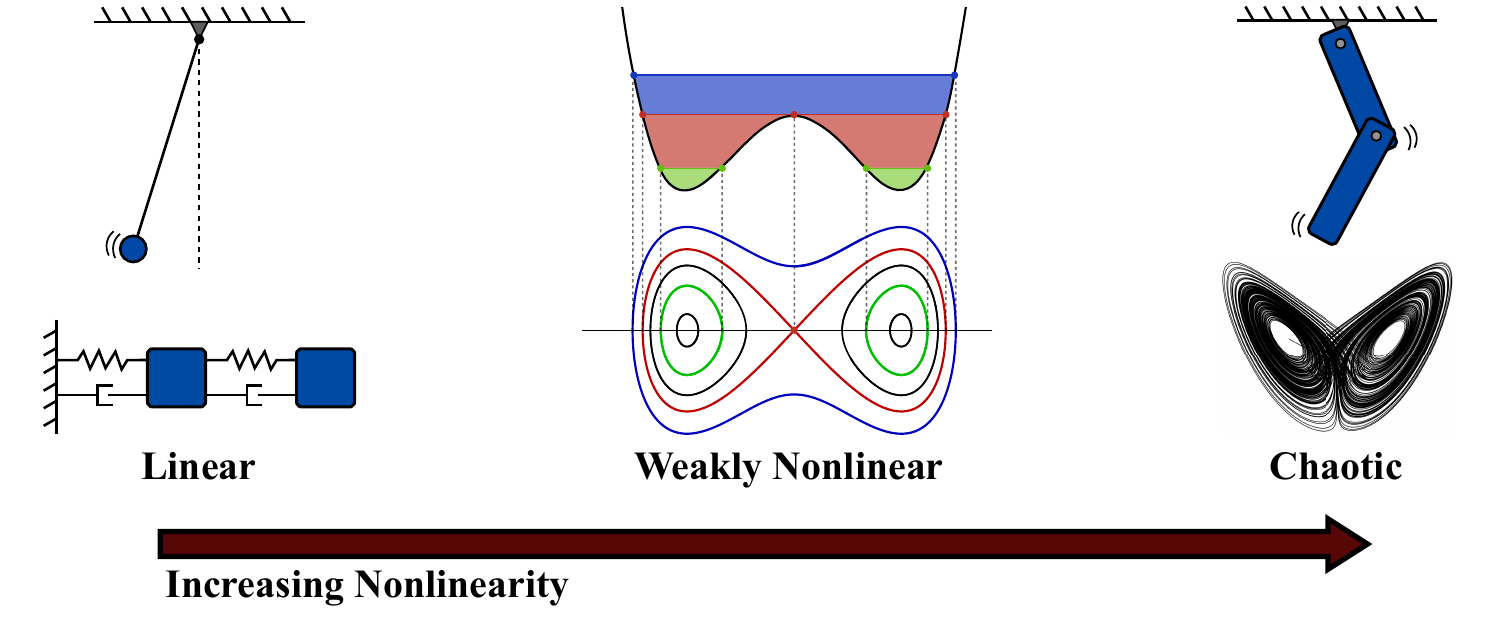}
\end{overpic}
\vspace{-.15in}
\caption{\small Chaotic dynamical systems are often viewed as a progression of increasing nonlinearity.}\label{fig03}
\end{center}
\vspace{-.25in}
\end{figure}

Historically, the two dominant perspectives on dynamical systems have either been geometric or statistical~\cite{Budivsic2012chaos}.  
In the geometric perspective, illustrated in Fig.~\ref{fig03}, the organization and topology of trajectories in phase space provides a qualitative picture of global dynamics and enables detailed quantitative descriptions of local dynamics near fixed points or periodic orbits~\cite{guckenheimer_holmes,Marsden1976book,MarsdenMTAA,Chorin1990book,Marsden:MS,koon2000heteroclinic}.  
Phase space transport is largely mediated by saddle points, and even in relatively simple systems such as the double pendulum or Lorenz system in Fig.~\ref{fig03}, the dynamics may give rise to chaotic dynamics.  
The statistical perspective trades the analysis of a single trajectory with the description of an ensemble of trajectories, providing a notion of mixing and uncertainty, while balancing the apparent structure and disorder in chaotic systems~\cite{Dellnitz2001book,Dellnitz2002hds,Dellnitz:2005,Froyland2009pd,Froyland2010chaos,Froyland2014siads,Kaiser2014jfm}.  
Recently, a third \emph{operator-theoretic} perspective, based on the evolution of measurement functions of the system, is gaining traction.  
This approach is not new, being introduced in 1931 by Koopman~\cite{Koopman1931pnas}, although the recent deluge of measurement data has renewed interest.  


Here, we develop a universal data-driven decomposition of chaos into a forced linear system. 
This relies on time-delay embedding, a cornerstone of dynamical systems, but takes a new perspective based on regression models~\cite{Brunton2016pnas} and modern Koopman operator theory~\cite{Mezic2005nd,Mezic2013arfm,Giannakis2015arxiv}.  
The resulting method partitions phase space into coherent regions where the forcing is small and dynamics are approximately linear, and regions where the forcing is large.  
The forcing may be measured from time series data and strongly predicts attractor switching and bursting phenomena in real-world examples. 
Linear representations of strongly nonlinear dynamics, enabled by machine learning and Koopman theory, promise to transform our ability to estimate, predict, and control complex systems in many diverse fields.  

\newpage
\section{Background}\label{sec:background}
The results in this paper are presented in the context of modern dynamical systems, specifically in terms of the Koopman operator.  
In this section, we provide a brief overview of relevant concepts in dynamical systems, including a discussion of Koopman operator theory in Sec.~\ref{Sec:Back:Koopman}, data-driven dynamical systems regression techniques in Sec.~\ref{Sec:Back:Regress}, and delay embedding theory in Sec.~\ref{Sec:Back:Delay}.

Throughout this work, we will consider dynamical systems of the form:
\begin{align}
\frac{d}{dt}\mathbf{x}(t) = \mathbf{f}(\mathbf{x}(t)).\label{Eq:ContinuousDynamics}
\end{align}
We will also consider the induced discrete-time dynamical system
\begin{eqnarray}
\mathbf{x}_{k+1}=\mathbf{F}(\mathbf{x}_k)\label{Eq:DiscreteDynamics}
\end{eqnarray}
where $\mathbf{x}_k$ may be obtained by sampling the trajectory in Eq.~\eqref{Eq:ContinuousDynamics} discretely in time, so that ${\mathbf{x}_k = \mathbf{x}(k\Delta t)}$.  The discrete-time propagator $\mathbf{F}$ is given by the flow map
\begin{eqnarray}
\mathbf{F}(\mathbf{x}_k)={\mathbf{x}_k+\int_{k\Delta t}^{({k+1})\Delta t}\mathbf{f}(\mathbf{x}(\tau))\,d\tau}.\label{Eq:FlowMap}
\end{eqnarray}
The discrete-time perspective is often more natural when considering experimental data.

\subsection{Koopman operator theory}\label{Sec:Back:Koopman}
Koopman spectral analysis was introduced in 1931 by B. O. Koopman~\cite{Koopman1931pnas} to describe the evolution of measurements of Hamiltonian systems, and this theory was generalized in 1932 by Koopman and von Neumann to systems with continuous spectra~\cite{Koopman1932pnas}. 
Koopman analysis provides an alternative to the more common geometric and statistical perspectives, instead describing the evolution operator that advances the space of measurement functions of the state of the dynamical system.  
The Koopman operator $\mathcal{K}$ is an infinite-dimensional linear operator that advances measurement functions $g$ of the state $\mathbf{x}$ forward in time according to the dynamics in \eqref{Eq:DiscreteDynamics}:
\begin{eqnarray}
\mathcal{K} g\triangleq g\circ\mathbf{F} \quad\Longrightarrow\quad \mathcal{K}g(\mathbf{x}_k) = g(\mathbf{x}_{k+1}).\label{Eq:Koopman}
\end{eqnarray}
Because this is true for \emph{all} measurement functions $g$, $\mathcal{K}$ is infinite dimensional and acts on the Hilbert space of functions of the state.  
For a detailed discussion on the Koopman operator, there are many excellent research articles~\cite{Mezic2005nd,Rowley2009jfm,Budivsic2009cdc,Budivsic2012physd,Lan2013physd} and reviews~\cite{Budivsic2012chaos,Mezic2013arfm}.

\begin{figure}
\begin{center}
\includegraphics[width=\textwidth]{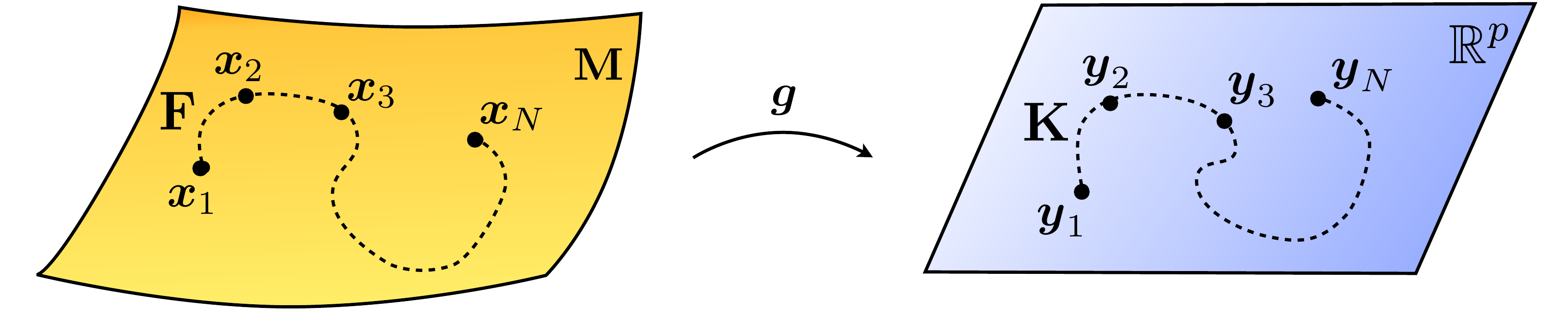}
\caption{Schematic illustrating the ability of the Koopman operator to globally linearize a nonlinear dynamical system with an appropriate choice of observable functions $g$.}\label{Fig:Koopman}
\end{center}
\end{figure}

A linear description of nonlinear dynamics is appealing, as many powerful analytic techniques exist to decompose, advance, and control linear systems.  
However, the Koopman framework trades finite-dimensional nonlinear dynamics for infinite-dimensional linear dynamics.  
Aside from a few notable exceptions~\cite{Gaspard1995pre,Bagheri2013jfm}, it is rare to obtain analytical representations of the Koopman operator.  
Obtaining a finite-dimensional approximation (i.e., a matrix $\mathbf{K}$) of the Koopman operator is therefore an important goal of data-driven analysis and control; this relies on a measurement subspace that remains invariant to the Koopman operator~\cite{Brunton2016plosone}.  Consider a measurement subspace spanned by measurement functions $\{g_1,g_2,\ldots,g_p\}$ so that for any measurement $g$ in this subspace 
\begin{eqnarray}
g = \alpha_1g_1 + \alpha_2 g_2 + \cdots + \alpha_p g_p\label{Eq:MeasurementSubspace}
\end{eqnarray}
then it remains in the subspace after being acted on by the Koopman operator
\begin{eqnarray}
\mathcal{K}g = \beta_1g_1 + \beta_2 g_2 + \cdots + \beta_p g_p.
\end{eqnarray}
In this case, we may restrict the Koopman operator to this $p$ dimensional measurement subspace and obtain a $p\times p$ matrix representation, $\mathbf{K}$, as illustrated in Fig.~\ref{Fig:Koopman}.  It has been shown previously that such a representation is useful for prediction and control of certain nonlinear systems that admit finite-dimensional Koopman invariant subspaces~\cite{Brunton2016plosone}.  If such a matrix representation exists, it is possible to define a linear system that advances the measurement functions, restricted to the subspace in \eqref{Eq:MeasurementSubspace}, as follows:
\begin{eqnarray}
\mathbf{y}_{k+1} = \mathbf{K}\mathbf{y}_{k},
\end{eqnarray}
where $\mathbf{y}_k=\begin{bmatrix}g_1(\mathbf{x}_k) & g_2(\mathbf{x}_k) & \cdots & g_p(\mathbf{x}_k\end{bmatrix}^T$ is a vector of measurements in the invariant subspace, evaluated at $\mathbf{x}_k$.  Left eigenvectors $\mathbf{\xi}$ of $\mathbf{K}$ give rise to Koopman eigenfunctions according to $\varphi = \mathbf{\xi}\mathbf{y}$.  

In practice, however, it is extremely challenging to obtain such a representation in terms of a Koopman invariant subspace.  
Moreover, it is impossible to obtain such an invariant subspace that contains linear measurements of the full state $\mathbf{x}$ for systems with more than one attractor, periodic orbit, and/or fixed point.  
This is simple to see, since a finite-dimensional linear system does not admit multiple fixed points or attracting structures.  
In addition, it is not always the case that the Koopman operator even has a discrete spectrum, as in mixing chaotic systems.  
However, the perspective of a data-driven linear approximation to a dynamical system is still valuable.  
Linear models can be obtained in entire basins of attraction of fixed points or periodic orbits using Koopman theory with the correct choice of measurement functions~\cite{Lan2013physd,Williams2015jnls}.  
Regression based methods to obtain a finite approximation of the Koopman operator, as described in Sec.~\ref{Sec:Back:Regress}, have become standard in the literature, although these methods all rely on a good choice of measurement functions.  
Identifying good measurement coordinates that approximately or exactly yield linear evolution equations will be one of the central challenges in dynamical systems in the coming years and decades.  
In the following, we will demonstrate the ability of delay coordinates to provide an approximately invariant measurement space for chaotic dynamics on an attractor.  

\subsection{Data-driven dynamic regression}\label{Sec:Back:Regress}
With increasingly large volumes of data, it is becoming possible to obtain models using modern regression techniques.  
This vibrant field will continue to grow as new techniques in machine learning make it possible to extract more information from data.  
In this section, we provide a brief overview of two leading regression-based system identification techniques:  1) the dynamic mode decomposition (DMD), which provides a best-fit linear operator from high-dimensional snapshot data, and may approximate the Koopman operator in some cases, and 2) the recent sparse identification of nonlinear dynamics (SINDy) algorithm, which produces parsimonious nonlinear models through sparse regression onto a library of nonlinear functions.  

\subsubsection{Dynamic mode decomposition (DMD)}\label{Sec:Back:Regress}
Dynamic mode decomposition (DMD) was originally introduced in the fluid dynamics community to decompose large experimental or numerical data sets into leading spatiotemporal coherent structures~\cite{Schmid2008aps,schmid:2010}.  
Shortly after, it was shown that the DMD algorithm provides a practical numerical framework to approximate the Koopman mode decomposition~\cite{Rowley2009jfm}.  
This connection between DMD and the Koopman operator was further strengthened and justified in a dynamic regression framework~\cite{Chen:2012,Tu2014jcd,Kutz2016book}.  

The DMD algorithm seeks a best-fit linear model to relate the following two data matrices
\begin{eqnarray}
 && {\bf X} = \begin{bmatrix}
\vline & \vline & & \vline \\
\bx_1 & \bx_2 & \cdots & \bx_{m-1}\\
\vline & \vline & & \vline
\end{bmatrix}  \qquad\qquad {\bf X}' = \begin{bmatrix}
\vline & \vline & & \vline \\
\bx_2 & \bx_3 & \cdots & \bx_{m}\\
\vline & \vline & & \vline
\end{bmatrix}.\label{eq:DataCollection}
\end{eqnarray}
The matrix $\mathbf{X}$ contains snapshots of the system state in time, and $\mathbf{X}'$ is a matrix of the same snapshots advanced a single step forward in time.  These matrices may be related by a best-fit linear operator $\mathbf{A}$ given by
\begin{eqnarray}
\mathbf{X}' = \mathbf{A}\mathbf{X} \quad\Longrightarrow\quad \mathbf{A} \approx \mathbf{X}'\mathbf{X}^{\dagger},
\end{eqnarray}
where $\mathbf{X}^{\dagger}$ is the pseudo-inverse, obtained via the singular value decomposition (SVD).  
The matrix $\mathbf{A}$ is a best-fit linear operator in the sense that it minimizes the Frobenius norm error $\|\mathbf{X}'-\mathbf{A}\mathbf{X}\|_F$.  

For systems of moderately large dimension, the operator $\mathbf{A}$ is intractably large, and so instead of obtaining $\mathbf{A}$ directly, we often seek the leading eigendecomposition of $\mathbf{A}$:
\begin{enumerate}
\item Take the SVD of $\mathbf{X}$:  
\begin{eqnarray}
\mathbf{X}=\mathbf{U}\mathbf{\Sigma}\mathbf{V}^*.
\end{eqnarray}
Here, $^*$ denotes complex conjugate transpose.  Often, only the first $r$ columns of $\mathbf{U}$ and $\mathbf{V}$ are required for a good approximation, $\mathbf{X}\approx \mathbf{\tilde{U}}\mathbf{\tilde{\Sigma}}\mathbf{\tilde{V}}$, where $\tilde{~}$ denotes a rank-$r$ truncation.
\item Obtain the $r\times r$ matrix $\mathbf{\tilde{A}}$ by projecting $\mathbf{A}$ onto $\mathbf{\tilde{U}}$:
\begin{eqnarray}
\mathbf{\tilde{A}} = \mathbf{\tilde{U}}^*\mathbf{A}\mathbf{\tilde{U}} = \mathbf{\tilde{U}}^*\mathbf{X}'\mathbf{\tilde{V}}\mathbf{\tilde{\Sigma}}^{-1}.
\end{eqnarray}
\item Compute the eigendecomposition of $\mathbf{\tilde{A}}$: 
\begin{eqnarray}
\mathbf{\tilde{A}} \mathbf{W} = \mathbf{W} \mathbf{\Lambda}.
\end{eqnarray}
The eigenvalues in $\mathbf{\Lambda}$ are eigenvalues of the full matrix $\mathbf{A}$.
\item Reconstruct full-dimensional eigenvectors of $\mathbf{A}$, given by the columns of $\mathbf{\Phi}$:
\begin{eqnarray}
\mathbf{\Phi} = \mathbf{X}'\mathbf{\tilde{V}}\mathbf{\tilde{\Sigma}}^{-1}\mathbf{W}.
\end{eqnarray}
\end{enumerate}

DMD, in its original formulation, is based on linear measurements of the state $\mathbf{x}$ of the system, such as velocity measurements from particle image velocimetry (PIV).  
This means that the measurement function $g$ is the identity map on the state.  
Linear measurements are not rich enough for many nonlinear dynamical systems, and so DMD has recently been extended to an augmented measurement vector including nonlinear functions of the state~\cite{Williams2015jnls}.  
However, choosing the correct nonlinear measurements that result in an approximately closed Koopman-invariant measurement system is still an open problem.  
Typically, measurement functions are either determined using information from the dynamical system (i.e., using quadratic nonlinearities for the Navier-Stokes equations), or by a brute-force search in a particular basis of Hilbert space (i.e., searching for polynomial functions or radial basis functions).  

\subsubsection{Sparse identification of nonlinear dynamics (SINDy)}
A recently developed technique, the sparse identification of nonlinear dynamics (SINDy) algorithm, identifies the nonlinear dynamics in Eq.~\eqref{Eq:ContinuousDynamics} from measurement data~\cite{Brunton2016pnas}.  
The SINDy algorithm uses sparse regression~\cite{Tibshirani1996lasso} in a nonlinear function space to determine the few active terms in the dynamics.  
Earlier related methods based on compressed sensing have been used to predict catastrophes in dynamical systems~\cite{Wang2011prl}.  
There are alternative methods that employ symbolic regression (i.e., genetic programming~\cite{koza1999genetic}) to identify dynamics~\cite{Bongard2007pnas,Schmidt2009science}. 
This work is part of a growing literature that is exploring the use of sparsity in dynamics~\cite{Ozolicnvs2013pnas,Schaeffer2013pnas,mackey2014compressive} and dynamical systems~\cite{Bai2014aiaa,Proctor2014epj,Brunton2014siads}.

The SINDy algorithm is an equation-free method~\cite{Kevrekidis2003cms} to identify a dynamical system \eqref{Eq:ContinuousDynamics} from data, much as in the DMD algorithm above.    
The basis of the SINDy algorithm is the observation that for many systems of interest, the function $\bf{f}$ only has a few active terms, making it sparse in the space of possible functions.  
Instead of performing a brute-force search for the active terms in the dynamics, sparse regression makes it possible to efficiently identify the few non-zero terms.

To determine the function $\bf{f}$ from data, we collect a time-history of the state $\bx(t)$ and the derivative $\dot\bx(t)$; note that $\dot\bx(t)$ may be approximated numerically from $\bx$.  The data is sampled at several times $t_1, t_2, \cdots, t_m$ 
and arranged into two large matrices:
\begin{subequations}
\begin{eqnarray}
\bX 
= \overset{\text{\normalsize state}}{\left.\overrightarrow{\begin{bmatrix}
x_1(t_1) & x_2(t_1) & \cdots & x_n(t_1)\\
x_1(t_2) & x_2(t_2) & \cdots & x_n(t_2)\\
\vdots & \vdots & \ddots & \vdots \\
x_1(t_m) & x_2(t_m) & \cdots & x_n(t_m)
\end{bmatrix}}\right\downarrow}\begin{rotate}{270}\hspace{-.125in}time~~\end{rotate}\label{Eq:DataMatrix}\qquad\qquad
\dot\bX 
= \begin{bmatrix}
\dot x_1(t_1) & \dot x_2(t_1) & \cdots & \dot x_n(t_1)\\
\dot x_1(t_2) & \dot x_2(t_2) & \cdots & \dot x_n(t_2)\\
\vdots & \vdots & \ddots & \vdots \\
\dot x_1(t_m) & \dot x_2(t_m) & \cdots & \dot x_n(t_m)
\end{bmatrix}.\\\nonumber
\end{eqnarray}
\end{subequations}
Next, we construct an augmented library $\bTheta(\bX)$ consisting of candidate nonlinear functions of the columns of $\bX$.  
For example, $\bTheta(\bX)$ may consist of constant, polynomial and trigonometric terms:
\begin{eqnarray}
\bTheta(\bX) = 
\begin{bmatrix} 
~~\vline&\vline & \vline & \vline & & \vline & \vline & \vline& \vline &  ~~ \\
~~\mathbf{1}&\bX & \bX^{P_2} & \bX^{P_3} & \cdots & \sin(\bX) & \cos(\bX) & \sin(2\bX) & \cos(2\bX) & \cdots ~~\\
~~\vline &\vline & \vline & \vline & & \vline &\vline &\vline & \vline &  ~~
\end{bmatrix}.\label{Eq:NonlinearLibrary}
\end{eqnarray}

Each column of $\bTheta(\bX)$ is a candidate function for the right hand side of Eq.~\eqref{Eq:ContinuousDynamics}.  
Since only a few of these nonlinearities are likely active in each row of $\bf{f}$, sparse regression is used to determine the sparse vectors of coefficients $\bXi = \begin{bmatrix}\bxi_1 &\bxi_2 & \cdots & \bxi_n\end{bmatrix}$ indicating which nonlinearities are active.
\begin{eqnarray}
\dot\bX = \bTheta(\bX) \bXi.\label{Eq:SparseRegression}
\end{eqnarray}

Once $\bXi$ has been determined, a model of each row of the governing equations may be constructed as follows:
\begin{eqnarray}
\dot{\bx}_k = {\bf f}_k(\bx) = \bTheta(\bx^T)\bxi_k.\label{Eq:sparseRow}
\end{eqnarray}
We may solve for $\bXi$ in Eq.~\eqref{Eq:SparseRegression} using sparse regression.  
In many cases, we may need to normalize the columns of $\bTheta(\bX)$ first to ensure that the restricted isometry property holds~\cite{Wang2011prl}; this is especially important when the entries in $\bX$ are small, since powers of $\bX$ will be minuscule.

Note that in the case the the library $\mathbf{\Theta}$ contains linear measurements of the state, the SINDy method reduces to a linear regression, closely related to the DMD above, but with transposed notation.  
The SINDy algorithm also generalizes naturally to discrete-time formulations.

\subsection{Time-delay embedding }\label{Sec:Back:Delay}
It has long been observed that choosing good measurements is critical to modeling, predicting, and controlling dynamical systems.  
The concept of \emph{observability} in a linear dynamical system provides conditions for when the full-state of a system may be estimated from a time-history of measurements of the system, providing a rigorous foundation for dynamic estimation, such as the Kalman filter~\cite{Kalman1960jfe,Ho1965aac,Welch1995book,dp:book,sp:book}.  
Although observability has been extended to nonlinear systems~\cite{Hermann1977ieeetac}, significantly fewer results hold in this more general context.  

The Takens embedding theorem~\cite{Takens1981lnm} provides a rigorous framework for analyzing the information content of measurements of a nonlinear dynamical system.  
It is possible to enrich a measurement, $x(t)$, with time-shifted copies of itself, $x(t-\tau)$, which are known as delay coordinates.  
Under certain conditions, the attractor of a dynamical system in delay coordinates is \emph{diffeomorphic} to the original attractor in the original state space.   
This is truly remarkable, as this theory states that in some cases, it may be possible to reconstruct the entire attractor of a turbulent fluid from a time series of a single point measurement.  
Similar \emph{differential} embeddings may be constructed by using derivatives of the measurement.  
Takens embedding theory has been related to nonlinear observability~\cite{Aeyels1981jco,Aguirre2005jpa}, providing a much needed connection between these two important fields.  

Delay embedding has been widely used to analyze and characterize chaotic systems~\cite{Farmer1987prl,Crutchfield1987cs,Sugihara1990nature,Rowlands1992physD,Abarbanel1993rmp,Sugihara2012science,Ye2015pnas}.  
The use of generalized delay coordinates are also used for linear system identification with the eigensystem realization algorithm (ERA)~\cite{ERA:1985} and in climate science with the singular spectrum analysis (SSA)~\cite{Broomhead1989prsla} and nonlinear Laplacian spectrum analysis (NLSA)~\cite{Giannakis2012pnas}.  
All of these methods are based on a singular value decomposition of a Hankel matrix, which is discussed below.  

\subsubsection{Hankel matrix analysis}\label{sec:ERA}
Both the eigensystem realization algorithm (ERA)~\cite{ERA:1985} and the singular spectrum analysis (SSA)~\cite{Broomhead1989prsla} are based on the construction of a Hankel matrix from a time series of measurement data.  
In the following, we will present the theory for a single scalar measurement, although this framework generalizes to multiple input, multiple output (MIMO) problems.  

The following Hankel matrix $\mathbf{H}$ is formed from a time series of a measurement $y(t)$:
\begin{eqnarray}
\mathbf{H} = \begin{bmatrix}y(t_1) & y(t_2) &  \cdots & y(t_{p}) \\ 
y(t_2) & y(t_3) &\cdots & y(t_{p+1}) \\
\vdots & \vdots  &\ddots & \vdots \\ 
y(t_q) & y(t_{q+1})& \cdots & y(t_{m})\end{bmatrix}.\label{Eq:Hankel}
\end{eqnarray}
Taking the singular value decomposition (SVD) of the Hankel matrix,
\begin{eqnarray}
\mathbf{H} = \mathbf{U\Sigma V}^T,
\end{eqnarray}
yields a hierarchical decomposition of the matrix into \emph{eigen} time series given by the columns of $\mathbf{U}$ and $\mathbf{V}$.  
These columns are ordered by their ability to express the variance in the columns and rows of the matrix $\mathbf{H}$, respectively.  
When the measurement $y(t)$ comes from the impulse response of an \emph{observable} linear system, then it is possible to use the SVD of the matrix $\mathbf{H}$ to reconstruct an accurate model of the full dynamics.  
This ERA procedure is widely used in system identification, and it has been recently connected to DMD~\cite{Tu2014jcd,Proctor2016siads}.  
In the following, we will generalize system identification using the Hankel matrix to nonlinear dynamical systems via the Koopman analysis.

\newpage
\section{Decomposing chaos: Hankel alternative view of Koopman (HAVOK)}  
Obtaining linear representations for strongly nonlinear systems has the potential to revolutionize our ability to predict and control these systems.  
In fact, the linearization of dynamics near fixed points or periodic orbits has long been employed for \emph{local} linear representation of the dynamics~\cite{guckenheimer_holmes}.  
The Koopman operator is appealing because it provides a \emph{global} linear representation, valid far away from fixed points and periodic orbits, although previous attempts to obtain finite-dimensional approximations of the Koopman operator have had limited success.  
Dynamic mode decomposition (DMD)~\cite{schmid:2010,Rowley2009jfm,Kutz2016book} seeks to approximate the Koopman operator with a best-fit linear model advancing spatial measurements from one time to the next.  
However, DMD is based on linear measurements, which are not rich enough for many nonlinear systems.  
Augmenting DMD with nonlinear measurements may enrich the model, but there is no guarantee that the resulting models will be closed under the Koopman operator~\cite{Brunton2016plosone}.  

Instead of advancing instantaneous measurements of the state of the system, we obtain intrinsic measurement coordinates based on the time-history of the system.  
This perspective is data-driven, relying on the wealth of information from previous measurements to inform the future.  
Unlike a linear or weakly nonlinear system, where trajectories may get trapped at fixed points or on periodic orbits, chaotic dynamics are particularly well-suited to this analysis: trajectories evolve to densely fill an attractor, so more data provides more information.

This method is shown in Fig.~\ref{Fig:Overview} for the Lorenz system in Sec.~\ref{Sec:Lorenz} below.  
The conditions of the Takens embedding theorem are satisfied~\cite{Takens1981lnm}, so eigen-time-delay coordinates may be obtained from a time series of a single measurement $x(t)$ by taking a singular value decomposition (SVD) of the following Hankel matrix $\mathbf{H}$:
\begin{eqnarray}
~\nonumber\\
\mathbf{H} = \begin{bmatrix}x(t_1) & x(t_2) &  \cdots & x(t_{p}) \\ 
x(t_2) & x(t_3) &\cdots & x(t_{p+1}) \\
\vdots & \vdots  &\ddots & \vdots \\ 
x(t_q) & x(t_{q+1})& \cdots & x(t_{m})
 \end{bmatrix} = \mathbf{U\Sigma V^*}.\label{Eq:Hankel}\\
 ~\nonumber
\end{eqnarray}
The columns of $\mathbf{U}$ and $\mathbf{V}$ from the SVD are arranged hierarchically by their ability to model the columns and rows of $\mathbf{H}$, respectively.  
Often, $\mathbf{H}$ may admit a low-rank approximation by the first $r$ columns of $\mathbf{U}$ and $\mathbf{V}$.  
Note that the Hankel matrix in \eqref{Eq:Hankel} is the basis of ERA~\cite{ERA:1985} in linear system identification and SSA~\cite{Broomhead1989prsla} in climate time series analysis.

The low-rank approximation to \eqref{Eq:Hankel} provides a \emph{data-driven} measurement system that is approximately invariant to the Koopman operator  for states on the attractor. 
By definition, the dynamics map the attractor into itself, making it \emph{invariant} to the flow.  
We may re-write \eqref{Eq:Hankel} with the Koompan operator $\mathcal{K}$:
\begin{eqnarray}
~\nonumber\\
\mathbf{H} = \begin{bmatrix}x(t_1) & \koop x(t_1) &  \cdots & \koop^{p-1}x(t_{1}) \\ 
\koop x(t_1) & \koop^2 x(t_1)  & \cdots & \koop^p x(t_{1}) \\
\vdots & \vdots  &\ddots & \vdots \\ 
\koop^{q-1} x(t_1) & \koop^q x(t_{1}) & \cdots & \koop^{m-1} x(t_{1})
 \end{bmatrix}.\label{Eq:HAVOK}\\
 ~\nonumber
\end{eqnarray}

\newpage
The columns of \eqref{Eq:Hankel}, and thus \eqref{Eq:HAVOK}, are well-approximated by the first $r$ columns of $\mathbf{U}$, so these eigen time series provide a Koopman-invariant measurement system.  
The first $r$ columns of $\mathbf{V}$ provide a time series of the magnitude of each of the columns of $\mathbf{U\Sigma}$ in the data.  
By plotting the first three columns of $\mathbf{V}$, we obtain an embedded attractor for the Lorenz system, shown in Fig.~\ref{Fig:Overview}.  
The rank $r$ can be obtained by the optimal hard threshold of Gavish and Donoho~\cite{Gavish2014ieeetit} or by other attractor dimension arguments~\cite{Abarbanel1993rmp}.  

The connection between eigen-time-delay coordinates from \eqref{Eq:Hankel} and the Koopman operator motivates a linear regression model on the variables in $\mathbf{V}$.  
Even with an approximately Koopman-invariant measurement system, there remain challenges to identifying a linear model for a chaotic system.  
A linear model, however detailed, cannot capture multiple fixed points or the unpredictable behavior characteristic of chaos with a positive Lyapunov exponent~\cite{Brunton2016plosone}.  
Instead of constructing a closed linear model for the first $r$ variables in $\mathbf{V}$, we build a linear model on the first $r-1$ variables and impose the last variable, $v_r$, as a forcing term.    
\begin{eqnarray}
\frac{d}{dt}\mathbf{v}(t) = \mathbf{A}\mathbf{v}(t) + \mathbf{B}{v}_r(t),\label{Eq:ChaosModel}
\end{eqnarray}
where $\mathbf{v}=\begin{bmatrix}v_1 & v_2 & \cdots & v_{r-1}\end{bmatrix}^T$ is a vector of the first ${r-1}$ eigen-time-delay coordinates.  
In all of the examples below, the linear model on the first $r-1$ terms is accurate, while no linear model represents $v_r$.  
Instead, $v_r$ is an input forcing to the linear dynamics in \eqref{Eq:ChaosModel}, which approximate the nonlinear dynamics in \eqref{Eq:ContinuousDynamics}.  
The statistics of $v_r(t)$ are non-Gaussian, as seen in the lower-right panel in Fig.~\ref{Fig:Overview}.  
The long tails correspond to rare-event forcing that drives lobe switching in the Lorenz system; this is related to rare-event forcing distributions observed and modeled by others~\cite{Majda2012nonlinearity,Sapsis2013pnas,Majda2014pnas}.  

The forced linear system in \eqref{Eq:ChaosModel} was discovered after applying the sparse identification of nonlinear dynamics (SINDy)~\cite{Brunton2016pnas} algorithm to delay coordinates of the Lorenz system.  
Even when allowing for the possibility of nonlinear dynamics for $\mathbf{v}$, the most parsimonious model was linear with a dominant off-diagonal structure in the $\mathbf{A}$ matrix (shown in Fig.~\ref{Fig:SparseModel}).   
This strongly suggests a connection with the Koopman operator, motivating the present work. 
The last term $v_r$ is not accurately represented by either linear or polynomial nonlinear models~\cite{Brunton2016pnas}.  
We refer to the framework presented here as the Hankel alternative view of Koopman (HAVOK) analysis.

The HAVOK analysis will be explored in detail below on the Lorenz system in Sec.~\ref{Sec:Lorenz} and on a wide range of numerical, experimental, and historical data models in Sec.~\ref{Sec:Results}.  
In nearly all of these examples, the forcing is generally small except for intermittent punctate events that correspond to transient attractor switching (for example, lobe switching in the Lorenz system) or bursting phenomena (in the case of Measles outbreaks).  
When the forcing signal is small, the dynamics are well-described by the Koopman linear system on the data-driven delay coordinates.  
When the forcing is large, the system is driven by an essential nonlinearity, which typically corresponds to an intermittent switching or bursting event.  
The regions of small and large forcing correspond to large coherent regions of phase space that may be analyzed further through machine learning techniques.

\newpage
\begin{figure}[b!]
\begin{center}
\begin{overpic}[width=\textwidth]{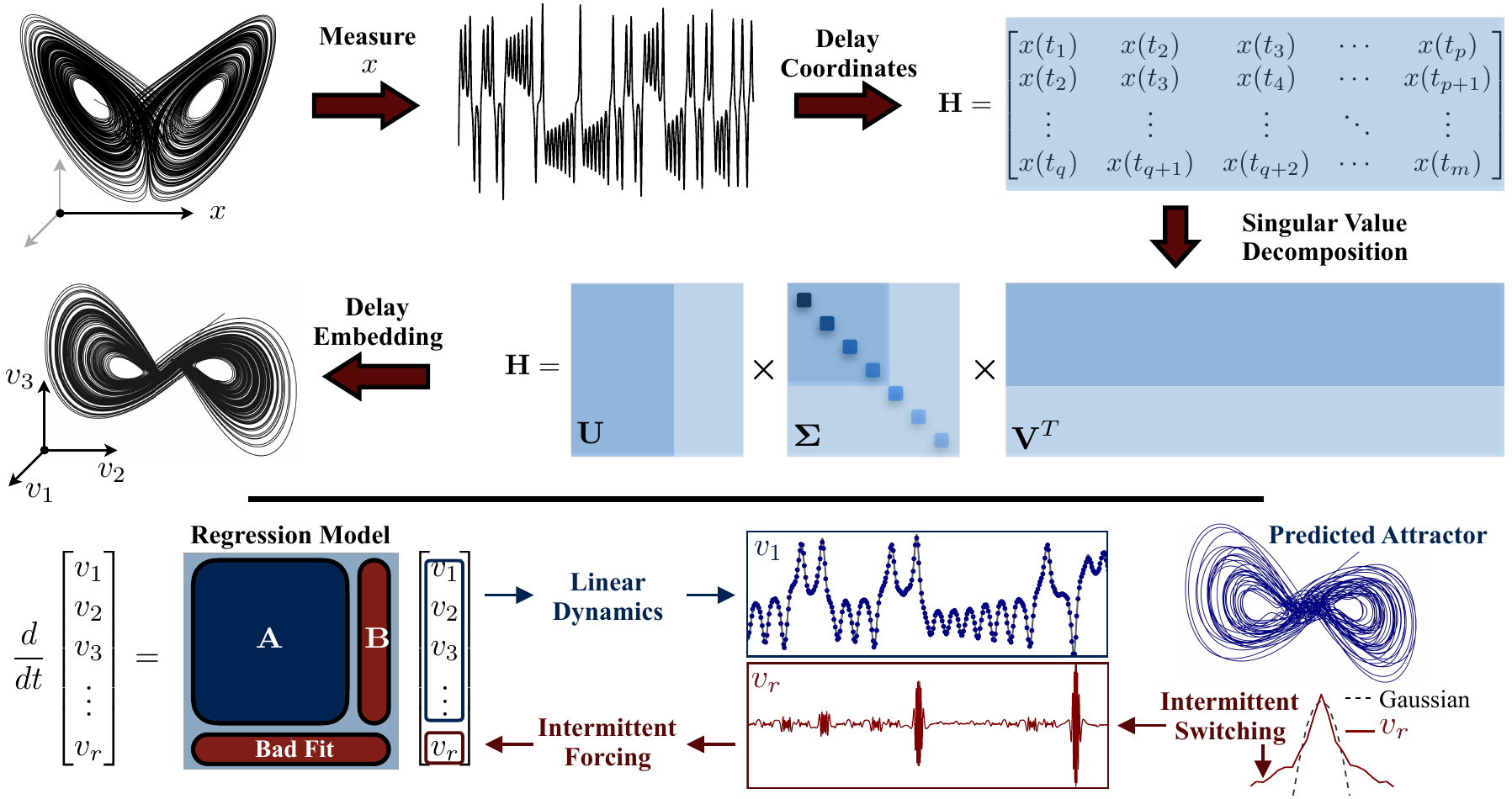}
\end{overpic}
\end{center}
\vspace{-.3in}
\caption{\small Decomposition of chaos into a linear dynamical system with intermittent forcing.  First, a time series $x(t)$ of of the chaotic Lorenz system is measured and time-shifted copies are stacked into a Hankel matrix $\mathbf{H}$.  Next, the singular value decomposition of $\mathbf{H}$ yields a hierarchical decomposition of  \emph{eigen} time series that characterize the measured dynamics.  In these coordinates, given by the columns of $\mathbf{V}$, we obtain a delay-embedded attractor.  Finally, a best-fit linear model is obtained on the time-delay coordinates $\mathbf{v}$; the linear fit for the first $r-1$ variables is excellent, but the last coordinate $v_r$ is not well-modeled with linear dynamics.  Instead, we consider $v_r$ as an input that forces the first $r-1$ variables.  The rare events in the forcing correspond to lobe switching in the chaotic dynamics.}\label{Fig:Overview}
\end{figure}
\section{HAVOK analysis illustrated on the chaotic Lorenz system}\label{Sec:Lorenz}
\begin{figure}[b]
\begin{center}
\begin{overpic}[width=.85\textwidth]{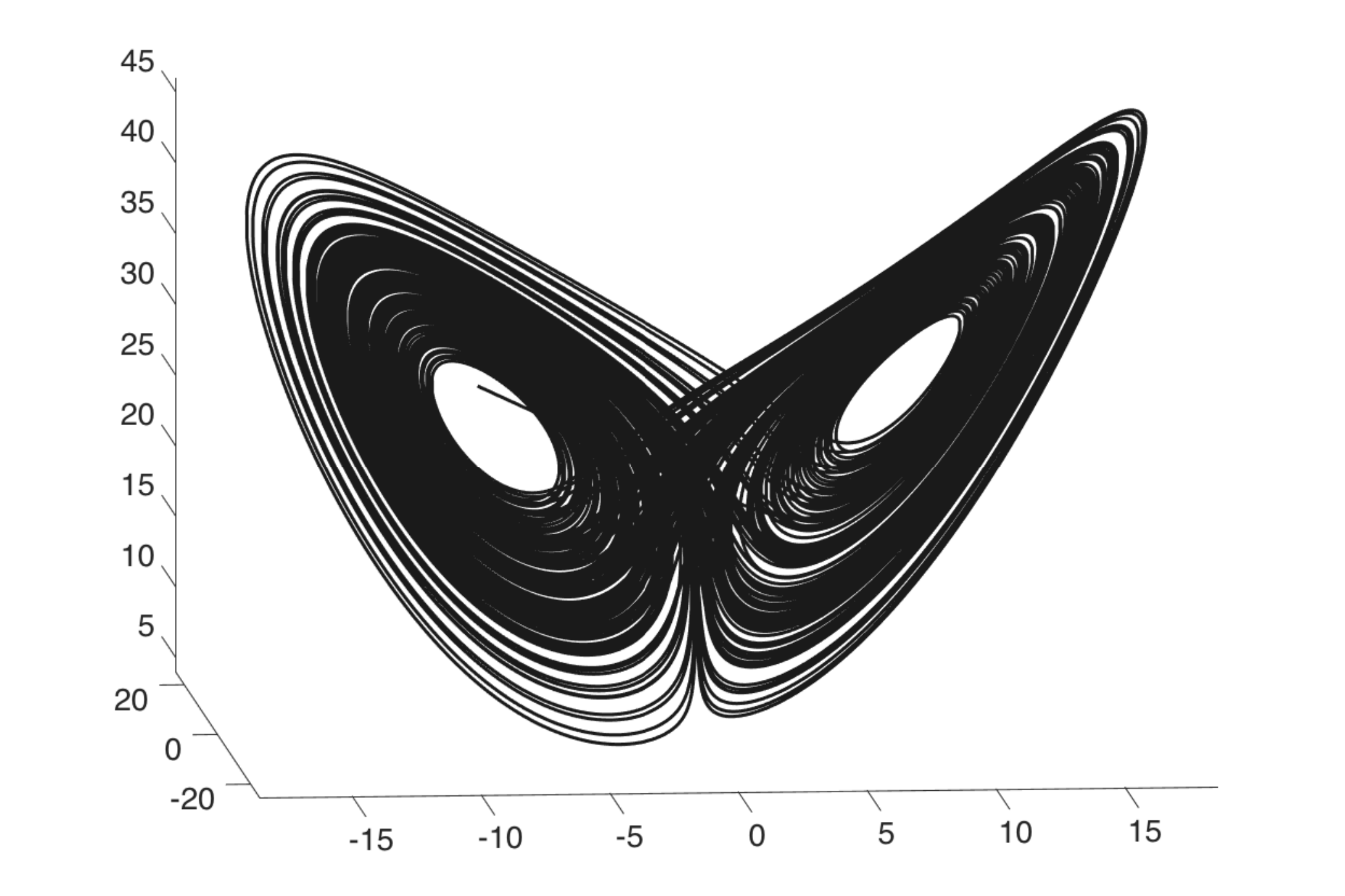}
\put(60,0){$x$}
\put(7,9){$y$}
\put(4,40){$z$}
\end{overpic}
\caption{Lorenz attractor from Eq.~\eqref{Eq:Lorenz}, simulated using the parameters in Tables.~\ref{Table:Parameters} and \ref{Table:HAVOK}.}\label{Fig:LorenzP1}
\end{center}
\end{figure}

To further understand the decomposition of a chaotic system into linear dynamics with intermittent forcing, we illustrate the HAVOK analysis in Fig.~\ref{Fig:Overview} on the chaotic Lorenz system~\cite{Lorenz1963jas}:
\begin{subequations}
\label{Eq:Lorenz}
\begin{eqnarray}
\dot{x} & = & \sigma (y - x)\\
\dot{y} & = & x(\rho -z) - y\\
\dot{z} & = & x y - \beta z,
\end{eqnarray}
\end{subequations}
with parameters $\sigma=10,\rho=28,$ and $\beta=8/3$.  
The Lorenz system is among the simplest and most well-studied examples of a deterministic dynamical system that exhibits chaos.  
A trajectory of the Lorenz system is shown in Fig.~\ref{Fig:LorenzP1}, integrated using the parameters in Tables~\ref{Table:Parameters} and \ref{Table:HAVOK} in Sec.~\ref{Sec:Results}.  
The trajectory moves along an attractor that is characterized by two lobes, switching back and forth between the two lobes intermittently by passing near a saddle point in the middle of the domain at $(x,y,z)=(0,0,0)$.  

\begin{figure}
\begin{center}
\begin{overpic}[width=\textwidth]{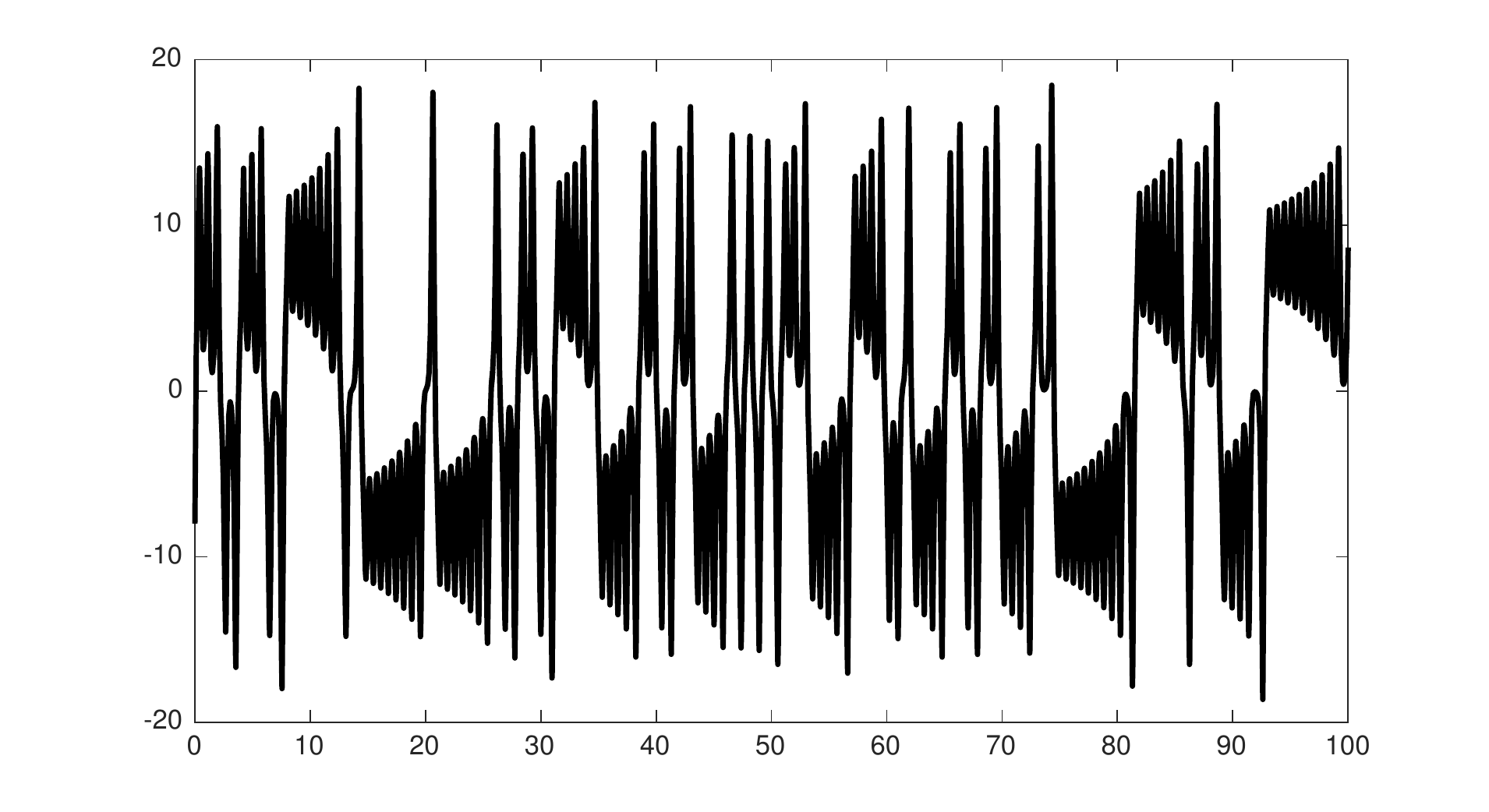}
\put(55,0){$t$}
\put(8,30){$x$}
\end{overpic}
\vspace{-.2in}
\caption{A time series of the variable $x(t)$ in the Lorenz system.  }\label{Fig:LorenzP2}
\vspace{-.2in}
\end{center}
\end{figure}

The various panels of Fig.~\ref{Fig:Overview} are provided in more detail below with labels and units.  
First, a time series of the $x$ variable is measured and plotted in Fig.~\ref{Fig:LorenzP2}.  
By stacking time-shifted copies of this measurement vector as rows of a Hankel matrix, as in Eq.~\eqref{Eq:Hankel}, it is possible to obtain an eigen-time-delay coordinate system through the singular value decomposition.  
These eigen-time-delay coordinates, given by the columns of the matrix $\mathbf{V}$, are the most self-similar time series features in the measurement $x(t)$, ordered by the variance they capture in the data.  

\begin{figure}[b]
\begin{center}
\begin{overpic}[width=.85\textwidth]{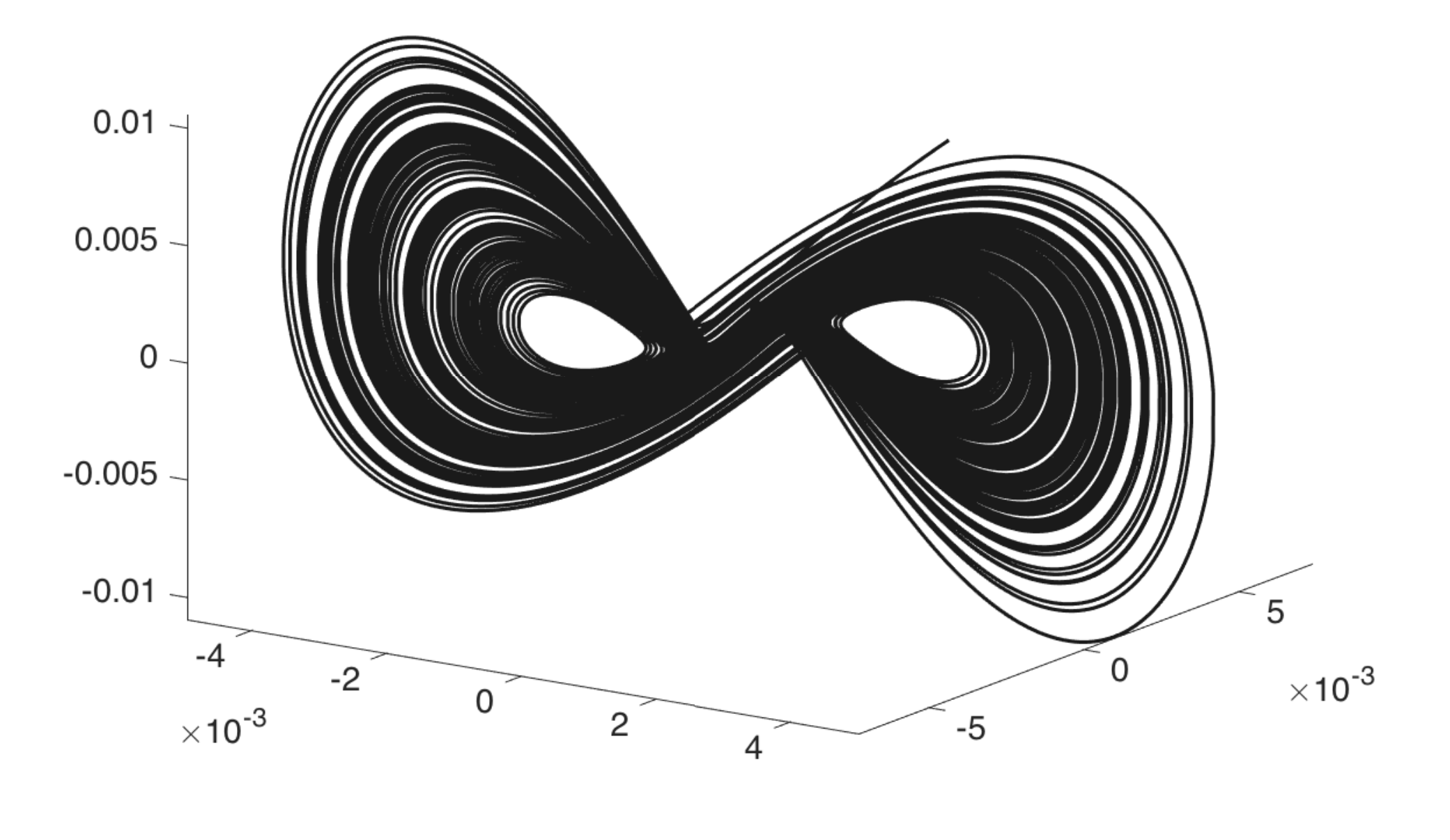}
\put(30,5){$v_1$}
\put(79,7){$v_2$}
\put(2,33){$v_3$}
\end{overpic}
\vspace{-.2in}
\caption{Time-delay embedded attractor using the eigen-time-delay coordinates obtained from the singular value decomposition of the Hankel matrix in Eq.~\eqref{Eq:Hankel}.}\label{Fig:LorenzP3}
\end{center}
\end{figure}

It is well known from Takens embedding theory~\cite{Takens1981lnm} that time delay coordinates provide an embedding of the scalar measurement $x(t)$ into a higher dimensional space.  The resulting attractor, shown in Fig.~\ref{Fig:LorenzP3}, is diffeomorphic to the original attractor in Fig.~\ref{Fig:LorenzP1}.

\begin{figure}[b]
\begin{center}
\begin{overpic}[width=.85\textwidth]{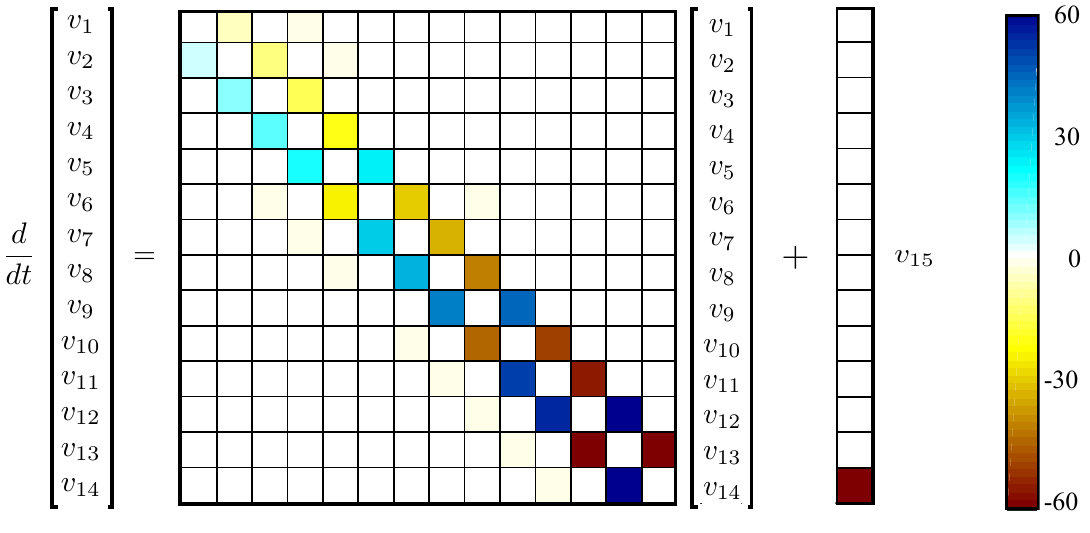}
\end{overpic}
\caption{HAVOK model obtained on time-delay coordinates of the Lorenz system.  Exact coefficients are provided in Tab.~\ref{Fig:HAVOKLorenzExact} in the Appendix.}\label{Fig:HAVOKModel}
\vspace{-.1in}
\end{center}
\end{figure}

Next, a HAVOK model is developed using the time delay coordinates.  In particular, a linear model is obtained for the first $14$ coordinates ($v_1,v_2,\dots,v_{14}$) with linear forcing from the $15^{\text{th}}$ coordinate $v_{15}$, given by the $15^{\text{th}}$ column of $\mathbf{V}$.  
This model is depicted schematically in Fig.~\ref{Fig:HAVOKModel}, and exact values are provided in Tab.~\ref{Fig:HAVOKLorenzExact} in the Appendix.  
The model may be obtained through a straightforward linear regression procedure, and an additional sparsity penalizing term may be added to eliminate terms in the model with very small coefficients~\cite{Brunton2016pnas}.  
The resulting model has a striking skew-symmetric structure, and the terms directly above and below the diagonal are nearly integer multiples of $5$.  
This fascinating structure is explored in Sec.~\ref{Sec:LorenzAdvanced}.

\begin{figure}[b]
\begin{center}
\begin{overpic}[width=.85\textwidth]{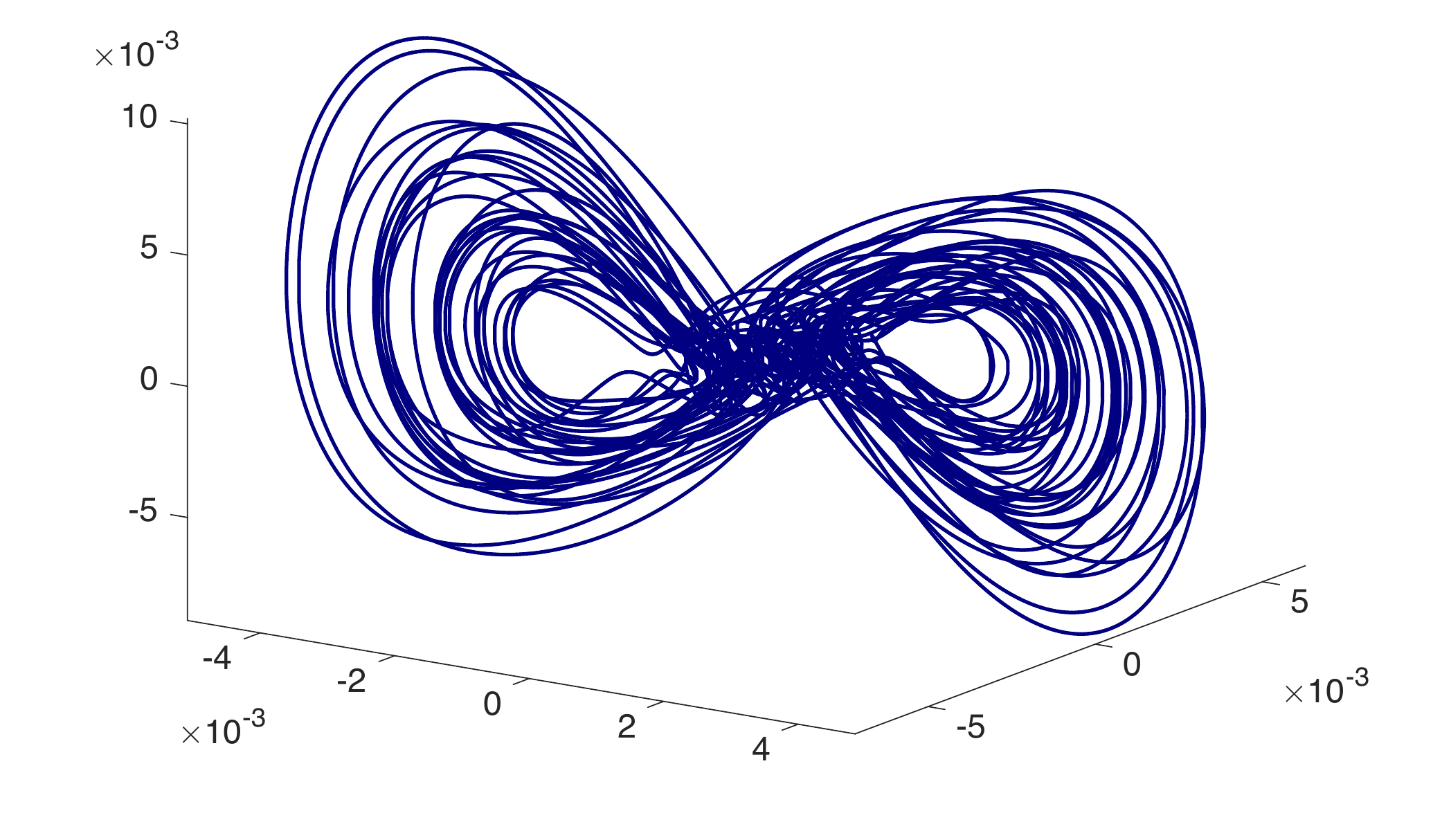}
\put(30,5){$v_1$}
\put(79,7){$v_2$}
\put(2,35){$v_3$}
\end{overpic}
\vspace{-.2in}
\caption{Reconstructed embedded attractor using linear HAVOK model with forcing from $v_{15}$.}\label{Fig:LorenzP5}
\end{center}
\end{figure}

\begin{figure}[b]
\begin{center}
\vspace{-.2in}
\begin{overpic}[width=\textwidth]{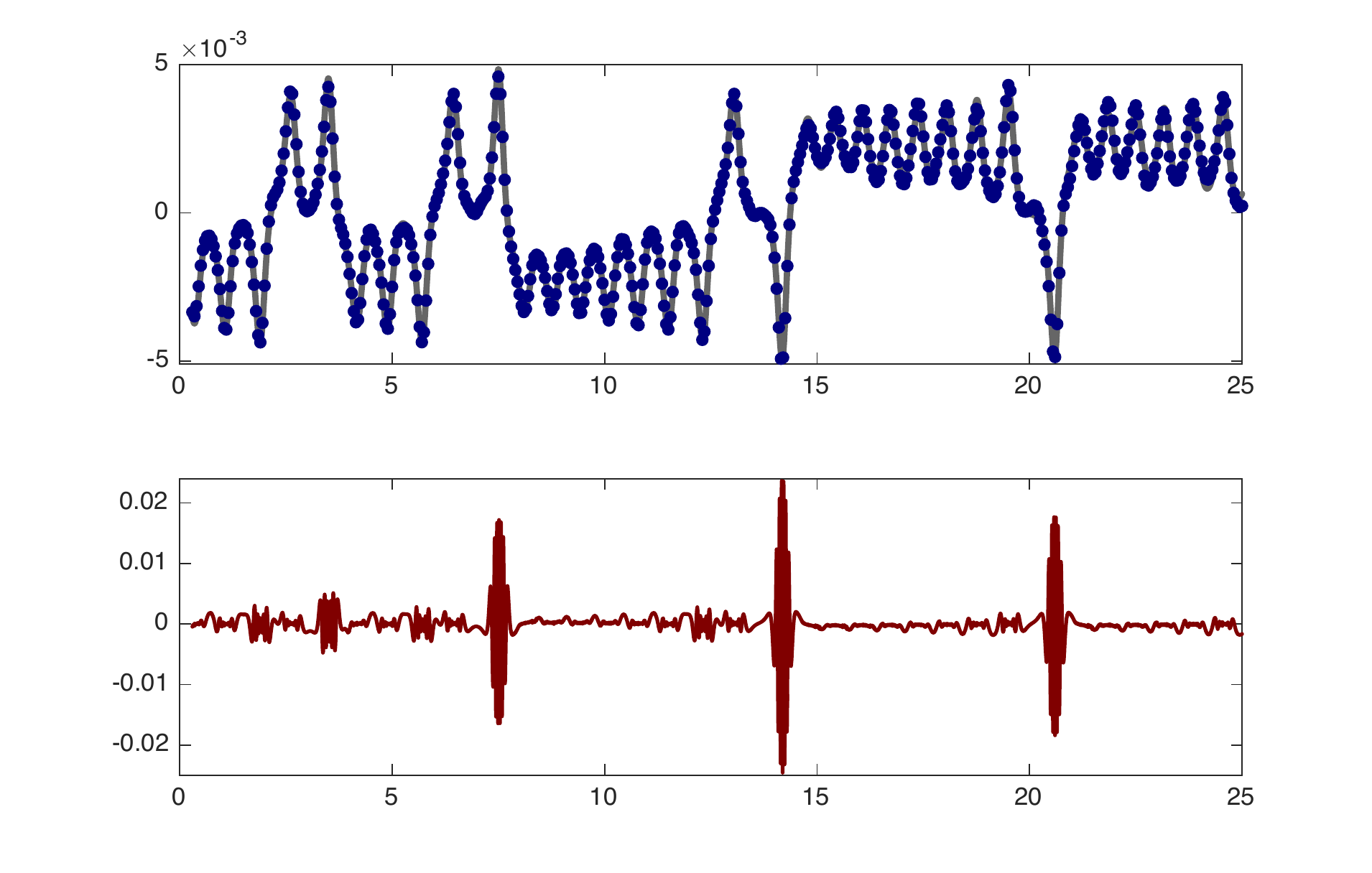}
\put(53,2){$t$}
\put(4,18){$v_{15}$}
\put(4,48){$v_{1}$}
\end{overpic}
\vspace{-.3in}
\caption{Reconstruction of $v_1$ using linear HAVOK model with forcing from $v_{15}$.}
\vspace{-.2in}
\label{Fig:LorenzP4}
\end{center}
\end{figure}

Using the HAVOK model with the signal $v_{15}$ as an input, it is possible to reconstruct the embedded attractor, as shown in Fig.~\ref{Fig:LorenzP5}.  
Figure~\ref{Fig:LorenzP4} shows the model prediction of the dominant time-delay coordinate, $v_1$, as well as the input forcing signal from $v_{15}$.  
From this figure, it is clear that for most times the forcing signal is nearly zero, and this signal begins to \emph{burst} when the trajectory is about to switch from one lobe of the attractor to another.  

\subsection{HAVOK models are predictive}
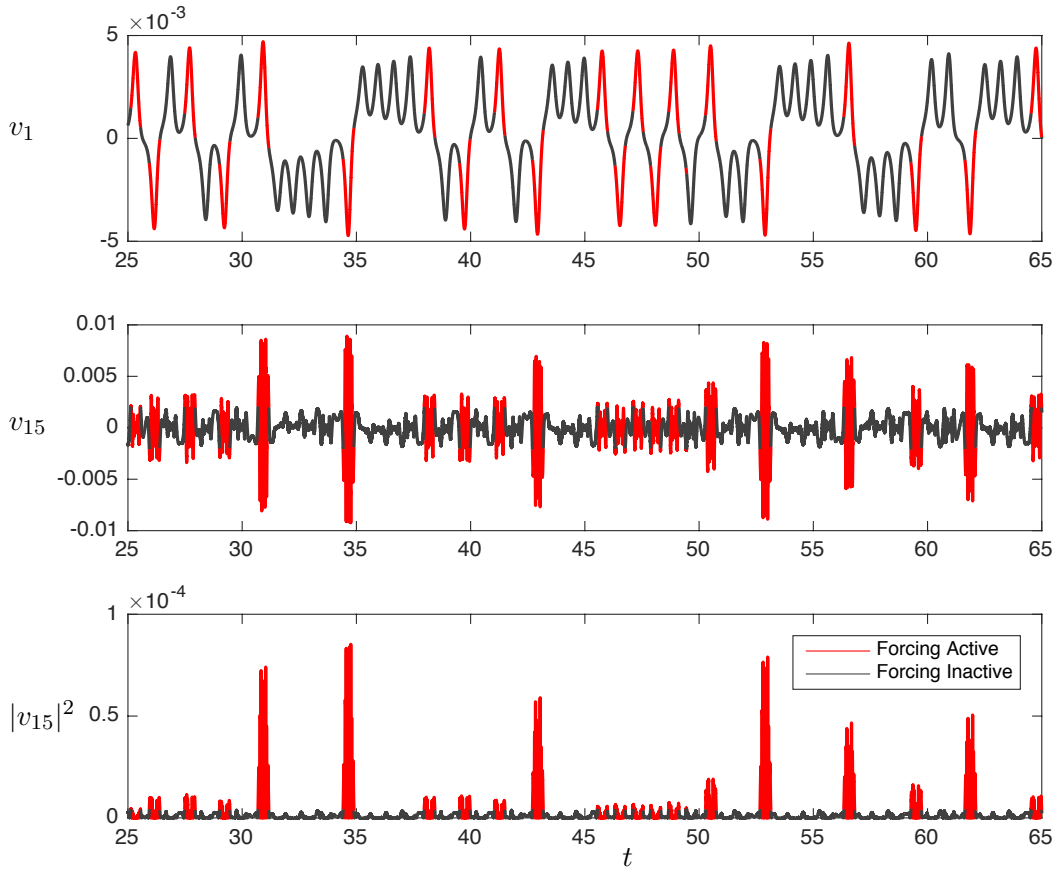
\begin{figure}[where!]
\begin{center}
\vspace{-.1in}
\begin{overpic}[width=.95\textwidth]{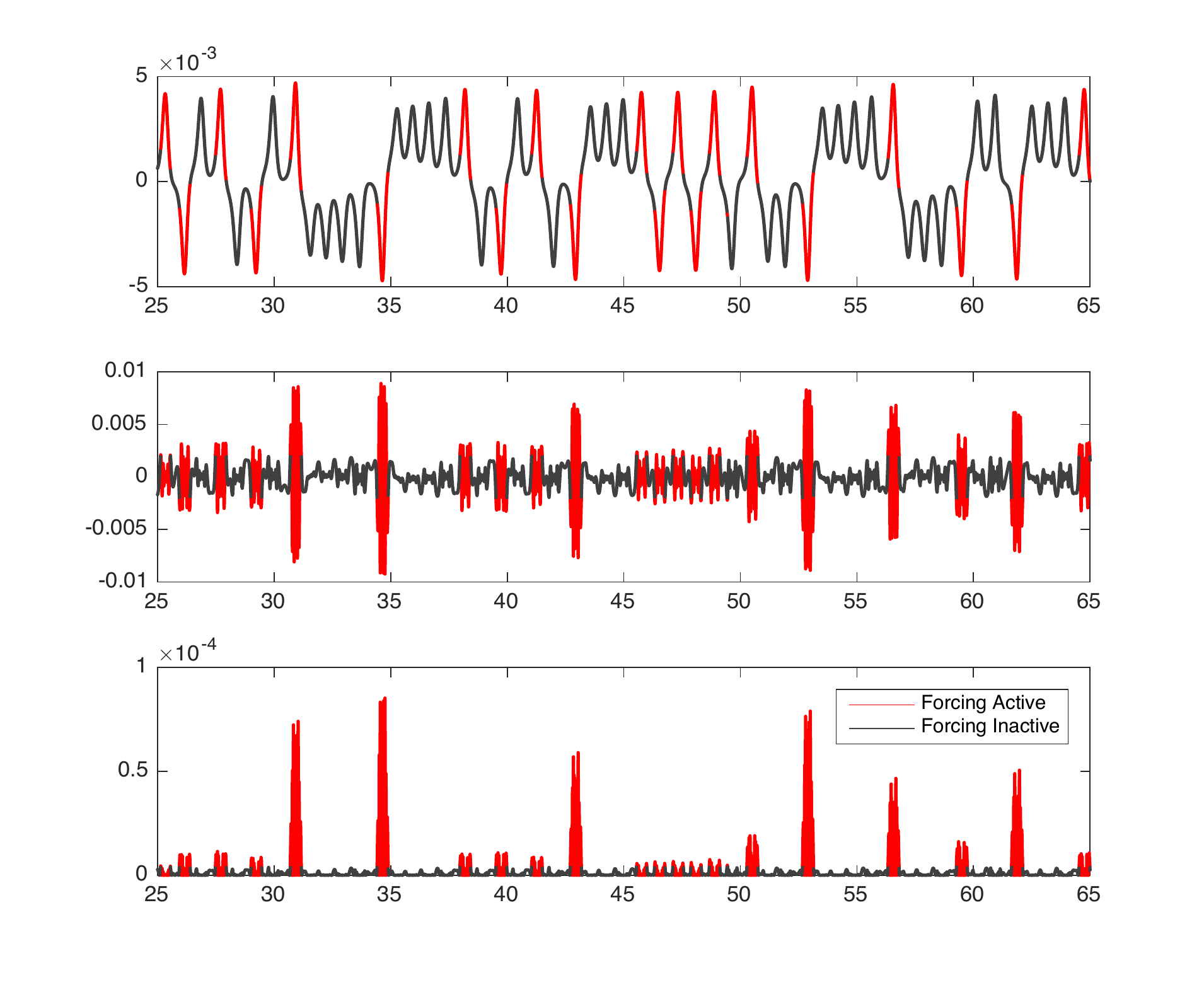}
\put(55,5){$t$}
\put(3,17){$|v_{15}|^2$}
\put(3,42){$v_{15}$}
\put(3,67){$v_1$}
\end{overpic}
\vspace{-.3in}
\caption{Delay coordinate $v_1$ of the Lorenz system, colored by the thresholded magnitude of the square of the forcing $v_{15}$.  When the forcing is active (red), the trajectory is about to switch, and when the forcing is inactive (grey), the solution is governed by predominantly linear dynamics.  }\label{Fig:LorenzP82D}
\vspace{-.15in}
\end{center}
\end{figure}

\begin{figure}
\begin{center}
\begin{overpic}[width=.85\textwidth]{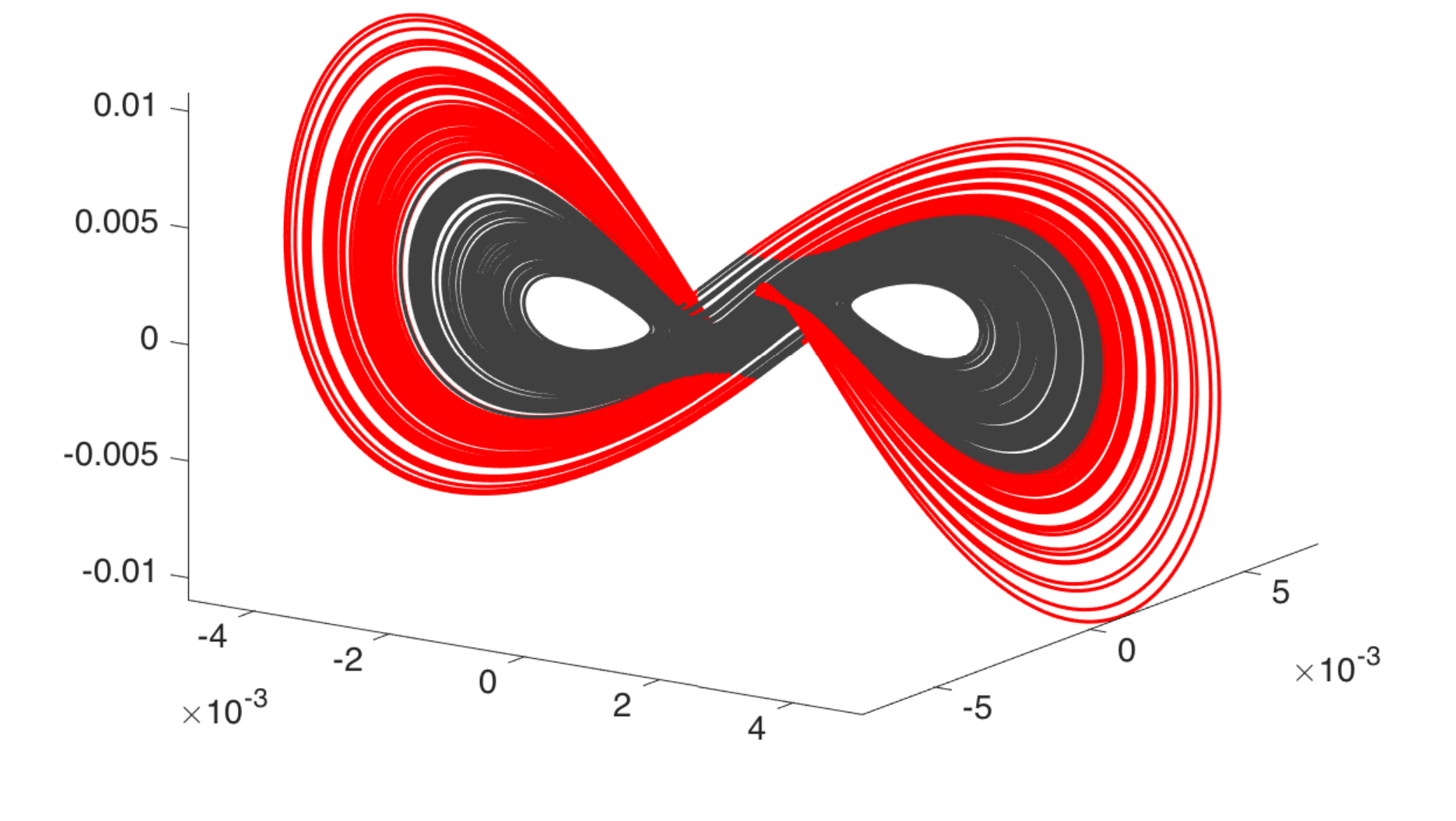}
\put(30,5){$v_1$}
\put(79,7){$v_2$}
\put(2,35){$v_3$}
\end{overpic}
\vspace{-.3in}
\caption{Time-delay embedded attractor of the Lorenz system color-coded by the activity of the forcing term $v_{15}$.  Trajectories in grey correspond to regions where the forcing is small and the dynamics are well approximated by Koopman linear dynamics.  The trajectories in red indicate that lobe switching is about to occur.}\label{Fig:LorenzP83D}
\end{center}
\end{figure}

\begin{figure}
\begin{center}
\vspace{-.1in}
\begin{overpic}[width=.9\textwidth]{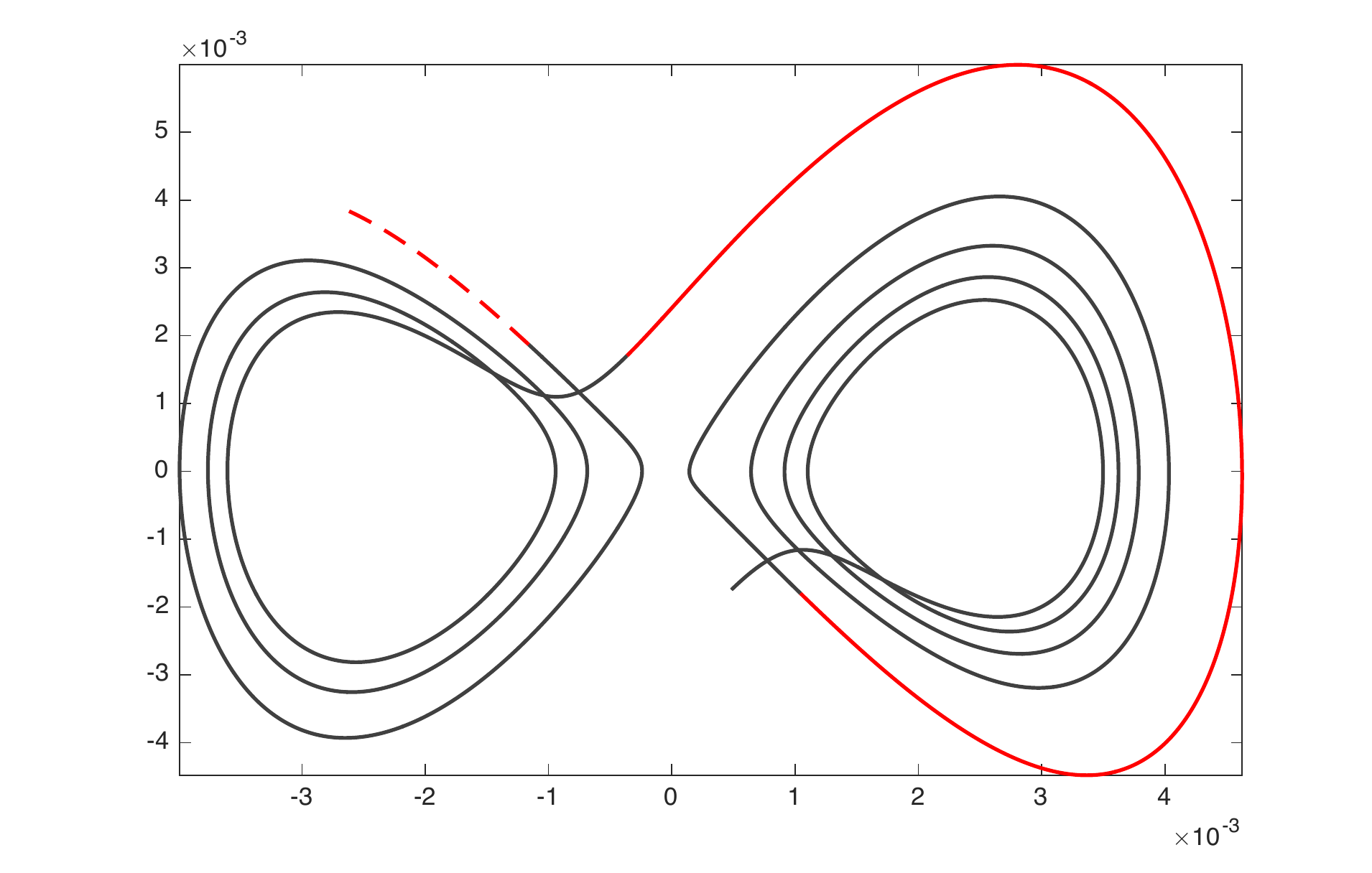}
\put(51,2){$v_1$}
\put(6,35){$v_2$}
\put(51,20){\textbf{A}}
\put(56.5,18.5){\textbf{B}}
\put(44.5,39.){\textbf{C}}
\put(38.,39.25){\textbf{D}}
\end{overpic}
\vspace{-.05in}
\caption{Illustration of one intermittent lobe switching event.  The trajectory starts at point (\textbf{A}), and resides in the basin of the right lobe for four revolutions, until (\textbf{B}), when the forcing becomes large, indicating an imminent switching event.  The trajectory makes one final revolution (red) and switches to the left lobe (\textbf{C}), where it makes three more revolutions.  At point (\textbf{D}), the activity of the forcing signal $v_{15}$ will increase, indicating that switching is imminent.}\label{Fig:LorenzP10}
\vspace{-.1in}
\end{center}
\end{figure}
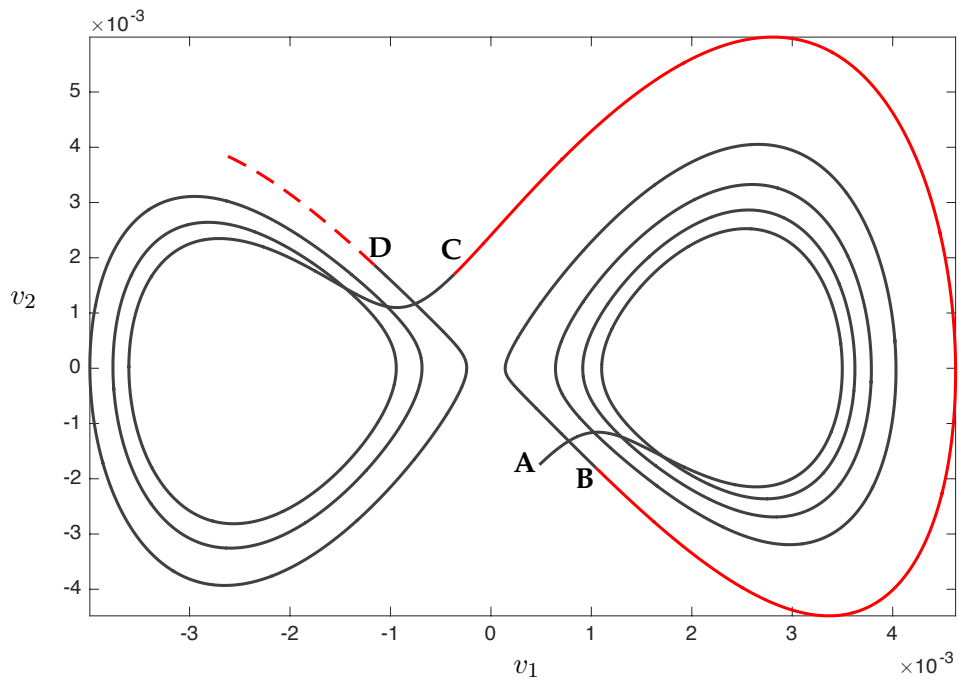

It is unclear from Fig.~\ref{Fig:LorenzP4} whether or not the bursting behavior in $v_{15}$ predicts the lobe switching in advance, or is simply active during the switch.  
To investigate this, it is possible to isolate regions where the forcing term $v_{15}$ is active by selecting values when $v_{15}^2$ are larger than a threshold value ($4\times 10^{-6}$), and coloring these sections of the trajectory in red.  The remaining portions of the trajectory, when the forcing is small, are colored in dark grey.  
These red and grey trajectory snippets are shown in Fig.~\ref{Fig:LorenzP82D}, where it is clear that the bursting significantly precedes the lobe switching.  
This is extremely promising, as it turns out that the activity of the forcing term is \emph{predictive}, preceding lobe switching by nearly one period. 
The same trajectories are plotted in three-dimensions in Fig.~\ref{Fig:LorenzP83D}, where it can be seen that the nonlinear forcing is active precisely when the trajectory is on the outer portion of the attractor lobes.
The prediction of lobe switching is quite good, although it can be seen that occasionally a switching event is missed by this procedure.  
A single lobe switching event is shown in Fig.~\ref{Fig:LorenzP10}, illustrating the geometry of the trajectories.  

It is important to note that when the forcing term is small, corresponding to the grey portions of the trajectory, the dynamics are largely governed by linear dynamics.  
Thus, the forcing term in effect distills the essential nonlinearity of the system, indicating when the dynamics are about to switch lobes of the attractor.  

\label{Sec:LorenzAdvanced}

\begin{figure}
\begin{center}
\vspace{-.2in}
\begin{overpic}[width=.9\textwidth]{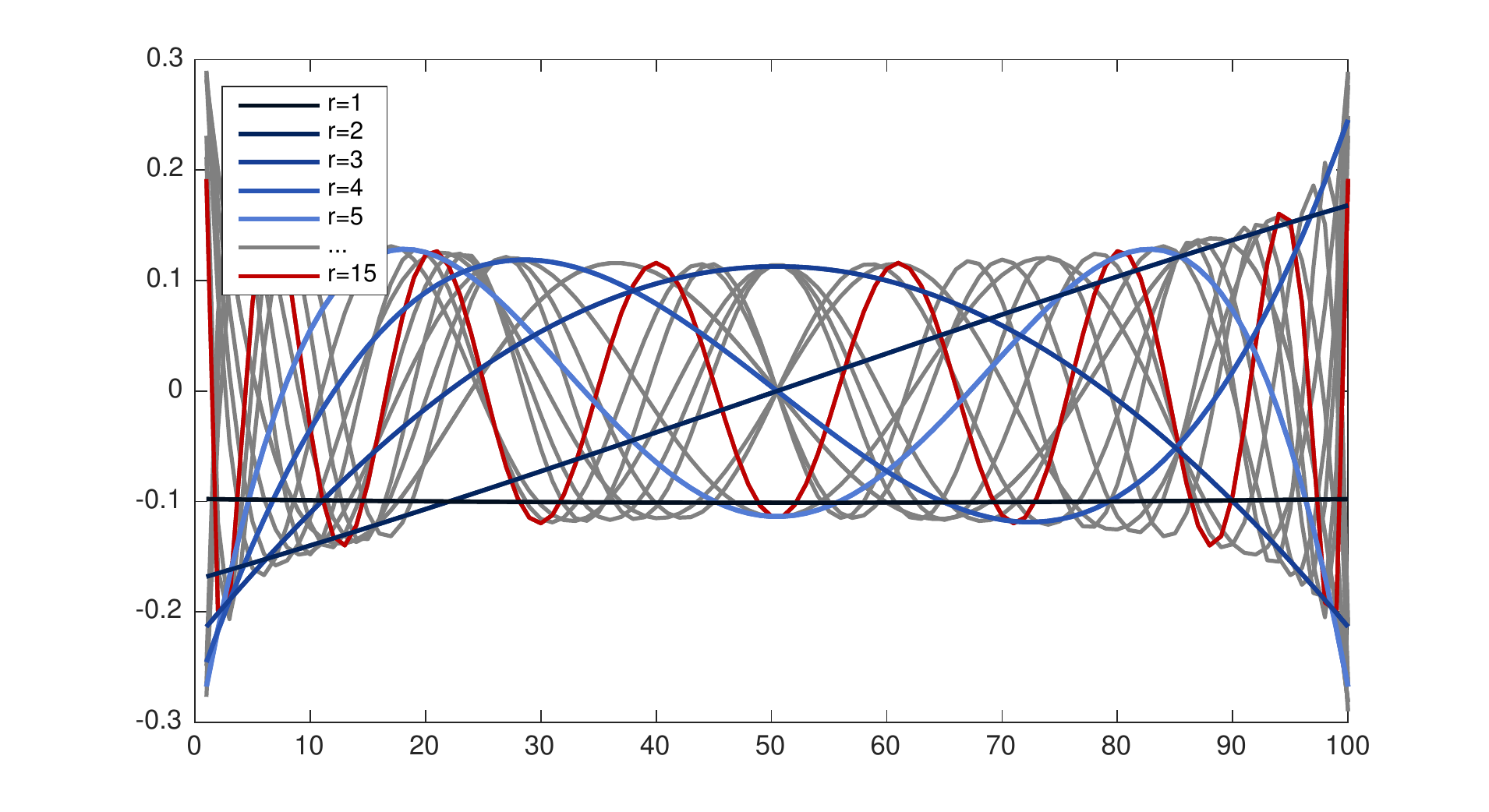}
\put(53,0){$t$}
\put(4,29){$u_r$}
\end{overpic}
\caption{Modes $u_r$ (columns of the matrix $\mathbf{U}$), indicating the short-time history that must be convolved with $x(t)$ to obtain $v_r$.}\label{Fig:LorenzP7}
\vspace{-.1in}
\end{center}
\end{figure}

\begin{figure}
\begin{center}
\begin{overpic}[width=.9\textwidth]{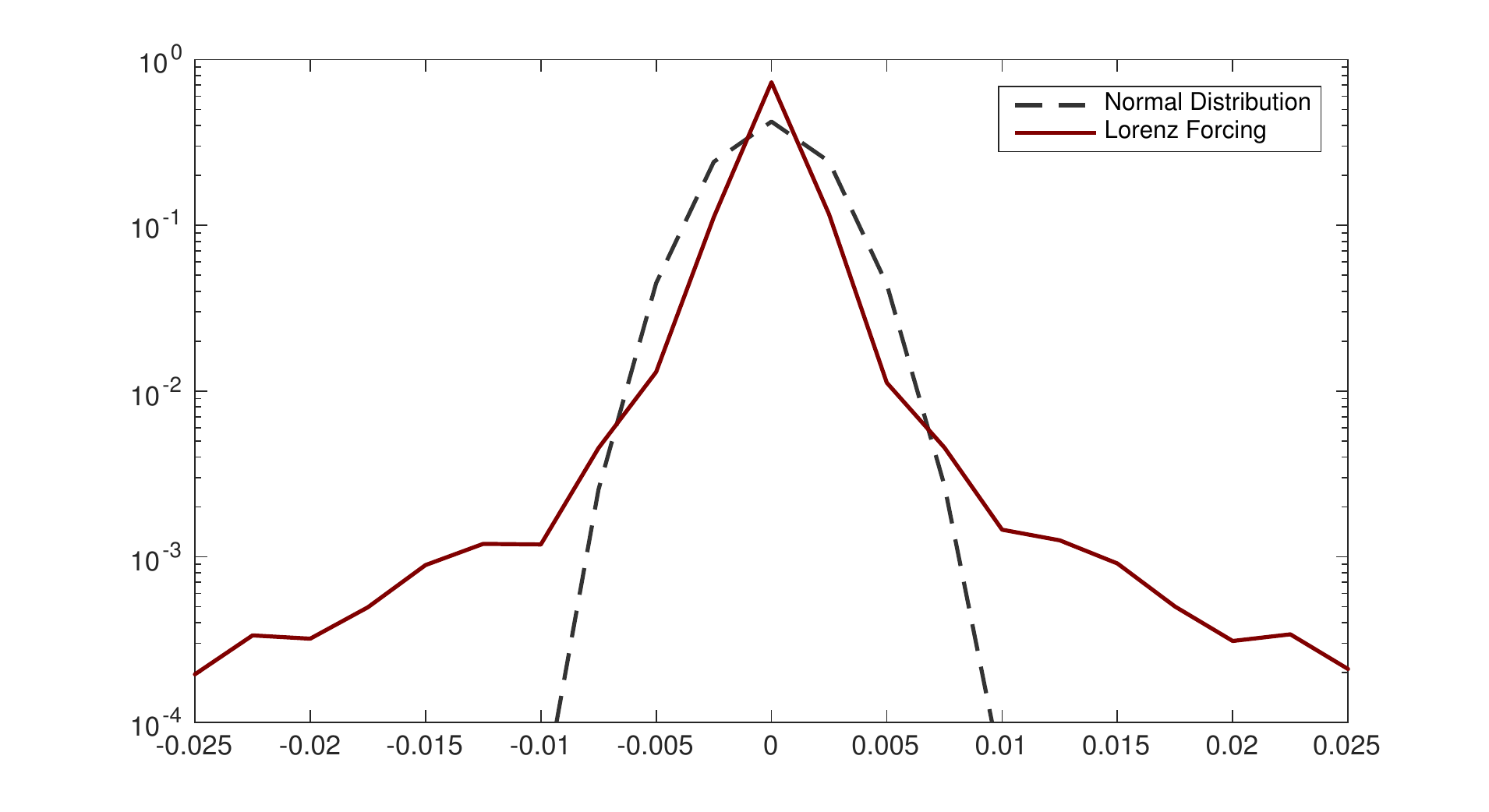}
\put(51,0){$v_{15}$}
\put(5,30){$p$}
\end{overpic}
\vspace{-.2in}
\caption{Probability density function of the forcing term in $v_{15}$ (red), plotted against the probability density function of an appropriately scaled normal distribution (black dash).}\label{Fig:LorenzP6}
\end{center}
\end{figure}

In practice, it is possible to measure $v_{15}$ from a streaming time series of $x(t)$ by convolution with the $u_{15}$ mode ($15^{\text{th}}$ column in the matrix $\mathbf{U}$), shown in Fig.~\ref{Fig:LorenzP7}.  
Finally, the probability density function of the forcing term $v_{15}$ is shown in Fig.~\ref{Fig:LorenzP6}; the long tails, compared with the normal distribution, indicate the rare intermittent switching events.   
Dynamic rare events are well-studied in climate and ocean waves~\cite{Cousins2014physd,Babaee2016prsa,Cousins2016jfm}.

\subsection{HAVOK models generalize beyond training data}
\begin{figure}[where]
\begin{center}
\begin{overpic}[width=\textwidth]{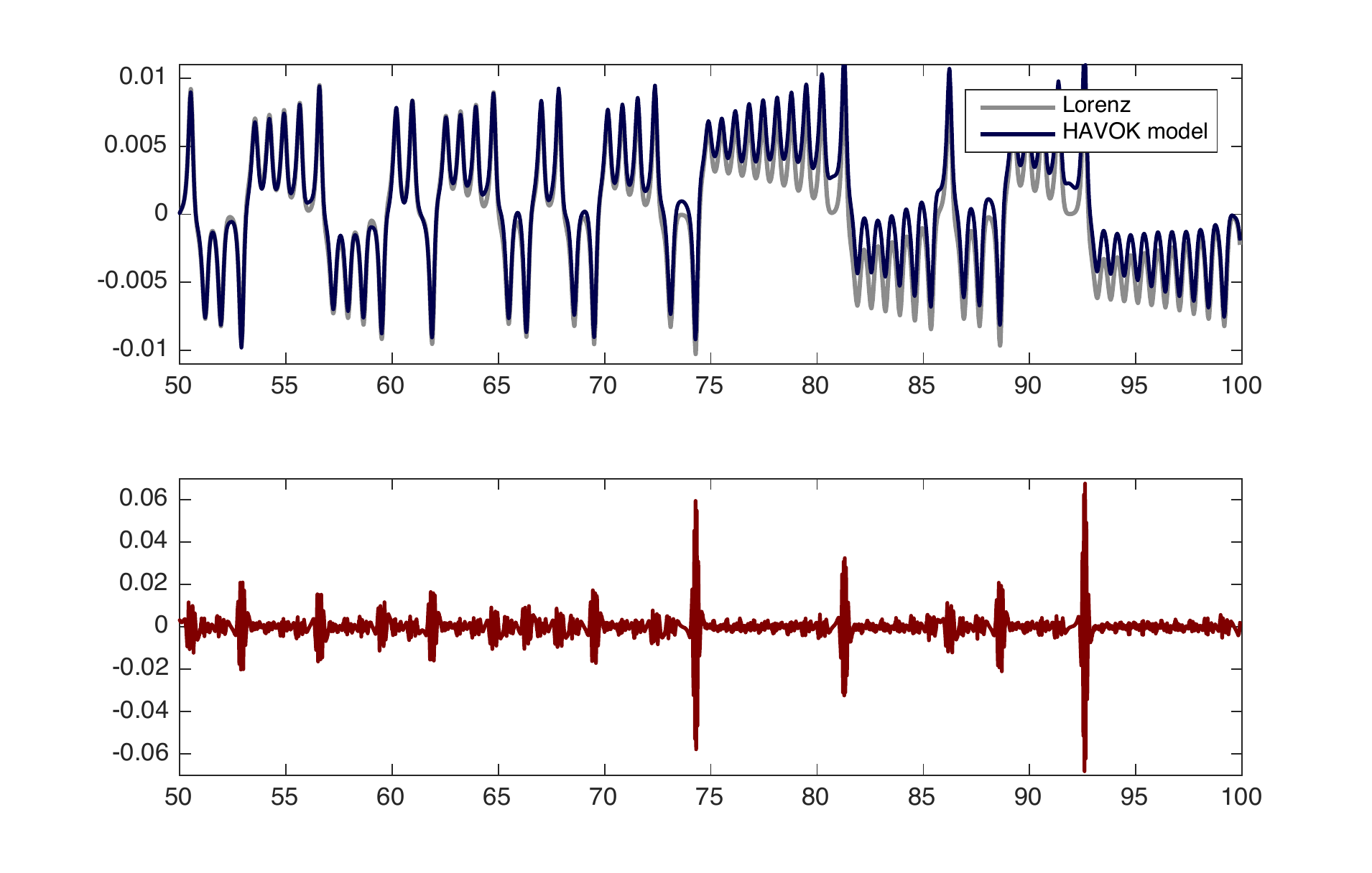}
\put(55,2.5){$t$}
\put(3,18){$v_{15}$}
\put(3,48){$v_{1}$}
\end{overpic}
\vspace{-.35in}
\caption{Performance of HAVOK model on test data from Lorenz system that was not used in the training set. This model was trained on a data set from $t=0$ to $t=50$.  }\label{Fig:LorenzP11}
\vspace{-.3in}
\end{center}
\end{figure}

It is important that the HAVOK model generalizes to test data that was not used to train the model.  Figure~\ref{Fig:LorenzP11} shows a validation test of a HAVOK model trained on data from $t=0$ to $t=50$ on new data from time $t=50$ to $t=100$.  
The model accurately captures the main lobe transitions, although small errors are gradually introduced for long time integration.  

\subsection{Structure and integrability of the HAVOK model for the Lorenz system}
The HAVOK model for the Lorenz system, provided in Tab.~\ref{Fig:HAVOKLorenzExact} in the Appendix, is highly structured.  
The $\mathbf{A}$ matrix of the linear dynamics is nearly skew-symmetric, and the terms directly above and below the diagonal are remarkably close to integer multiples of $5$.  
Intrigued by this structure, we have constructed an approximate system below that represents the \emph{idealized} structure in the HAVOK model for the Lorenz system:
\setcounter{MaxMatrixCols}{20}
\begin{eqnarray}
\scriptsize
\frac{d}{dt}
\begin{bmatrix}
v_1\\ v_2 \\ v_3 \\ v_4 \\ v_5 \\ v_6 \\ v_7 \\ v_8 \\ v_9 \\ v_{10} \\ v_{11}\\ v_{12}\\ v_{13} \\ v_{14}
\end{bmatrix} 
\hspace{-.05in}
=
\hspace{-.05in} 
\begin{bmatrix}
0 & -5 & 0 & 0 & 0 & 0 & 0 & 0 & 0 & 0 & 0 & 0 & 0 & 0\\
5 & 0 & -10 & 0 & 0 & 0 & 0 & 0 & 0 & 0 & 0 & 0 & 0 & 0\\
0 & 10 & 0 & -15 & 0 & 0 & 0 & 0 & 0 & 0 & 0 & 0 & 0 & 0\\
0 & 0 & 15 & 0 & -20 & 0 & 0 & 0 & 0 & 0 & 0 & 0 & 0 & 0\\
0 & 0 & 0 & 20 & 0 & 25 & 0 & 0 & 0 & 0 & 0 & 0 & 0 & 0\\
0 & 0 & 0 & 0 & -25 & 0 & -30 & 0 & 0 & 0 & 0 & 0 & 0 & 0\\
0 & 0 & 0 & 0 & 0 & 30 & 0 & -35 & 0 & 0 & 0 & 0 & 0 & 0\\
0 & 0 & 0 & 0 & 0 & 0 & 35 & 0 & -40 & 0 & 0 & 0 & 0 & 0\\
0 & 0 & 0 & 0 & 0 & 0 & 0 & 40 & 0 & 45 & 0 & 0 & 0 & 0\\
0 & 0 & 0 & 0 & 0 & 0 & 0 & 0 & -45 & 0 & -50 & 0 & 0 & 0\\
0 & 0 & 0 & 0 & 0 & 0 & 0 & 0 & 0 & 50 & 0 & -55 & 0 & 0\\
0 & 0 & 0 & 0 & 0 & 0 & 0 & 0 & 0 & 0 & 55 & 0 & 60 & 0\\
0 & 0 & 0 & 0 & 0 & 0 & 0 & 0 & 0 & 0 & 0 & -60 & 0 & -65\\
0 & 0 & 0 & 0 & 0 & 0 & 0 & 0 & 0 & 0 & 0 & 0 & 65 & 0
\end{bmatrix} 
\hspace{-.05in}
\begin{bmatrix}
v_1\\ v_2 \\ v_3 \\ v_4 \\ v_5 \\ v_6 \\ v_7 \\ v_8 \\ v_9 \\ v_{10} \\ v_{11}\\ v_{12}\\ v_{13} \\ v_{14}
\end{bmatrix}
\hspace{-.05in} 
+ 
\hspace{-.05in}
\begin{bmatrix}
0\\ 0\\0\\0\\0\\0\\0\\0\\0\\0\\0\\0\\0\\-70
\end{bmatrix}
v_{15}.
\end{eqnarray}

\begin{figure}[b]
\begin{center}
\vspace{-1in}
\begin{overpic}[width=\textwidth]{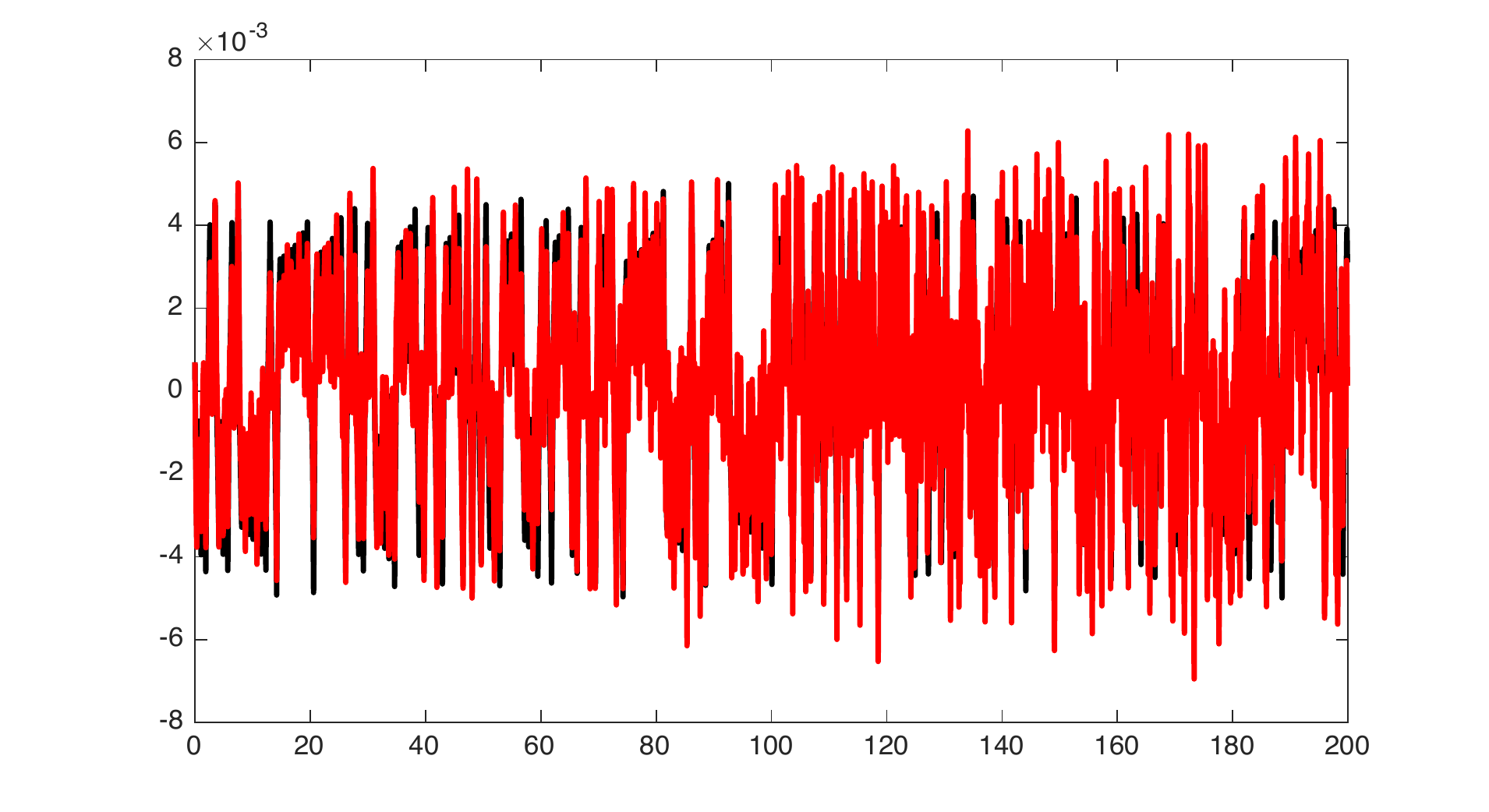}
\put(52,0){$t$}
\put(5.5,30){$v_1$}
\end{overpic}
\caption{Response of idealized HAVOK model with integer coefficients, forced by the $v_{15}$ time series from the Lorenz system.}\label{Fig:LorenzP9}
\vspace{-.2in}
\end{center}
\end{figure}

\begin{figure}[b]
\begin{center}
\begin{overpic}[width=\textwidth]{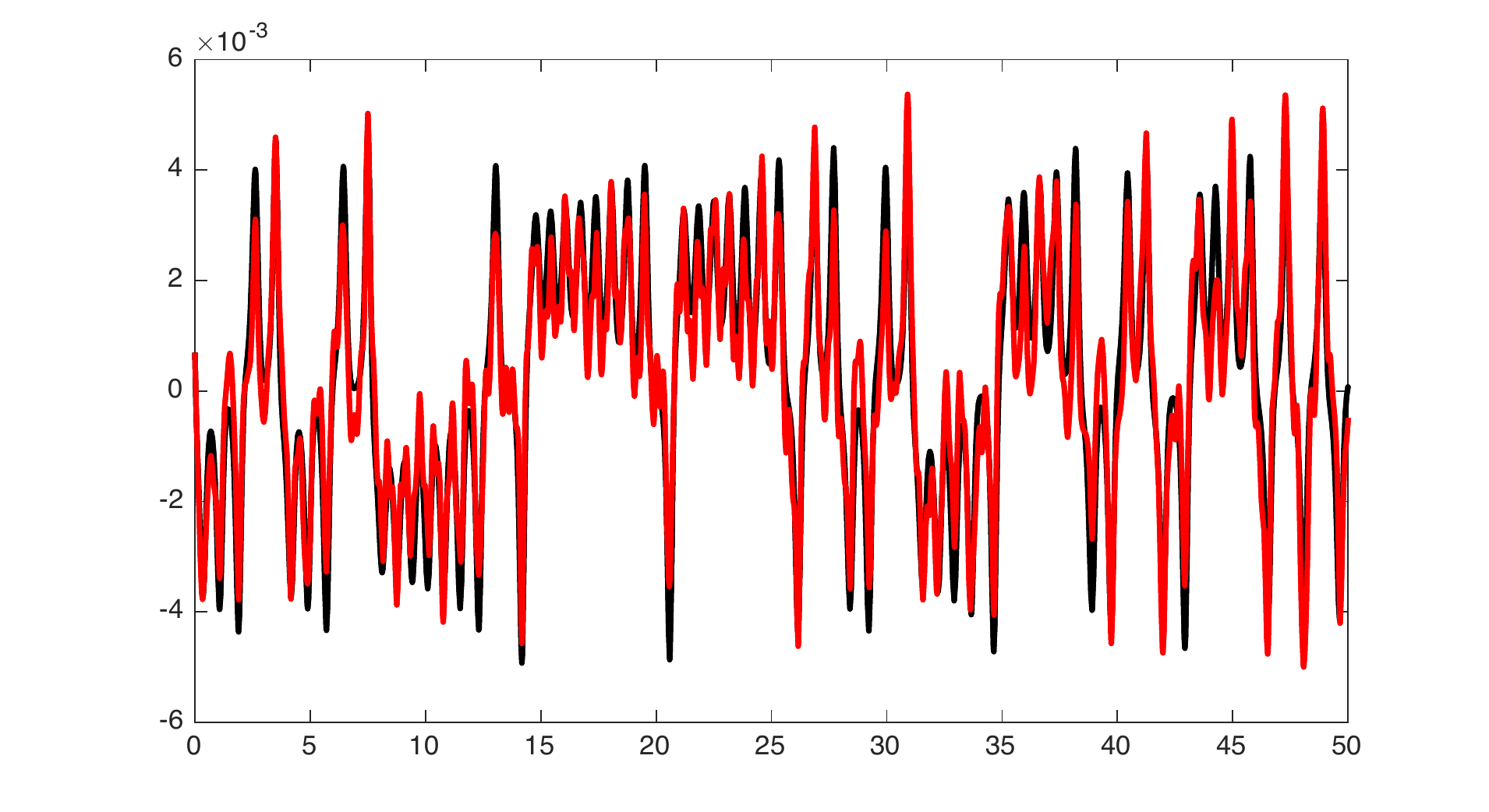}
\put(52,0){$t$}
\put(5.5,30){$v_1$}
\end{overpic}
\vspace{-.1in}
\caption{First 50 time unites of response of idealized HAVOK model with integer coefficients, forced by the $v_{15}$ time series from the Lorenz system.}\label{Fig:LorenzP9A}
\end{center}
\vspace{-.2in}
\end{figure}

\begin{figure}[t]
\begin{center}
\begin{overpic}[width=\textwidth]{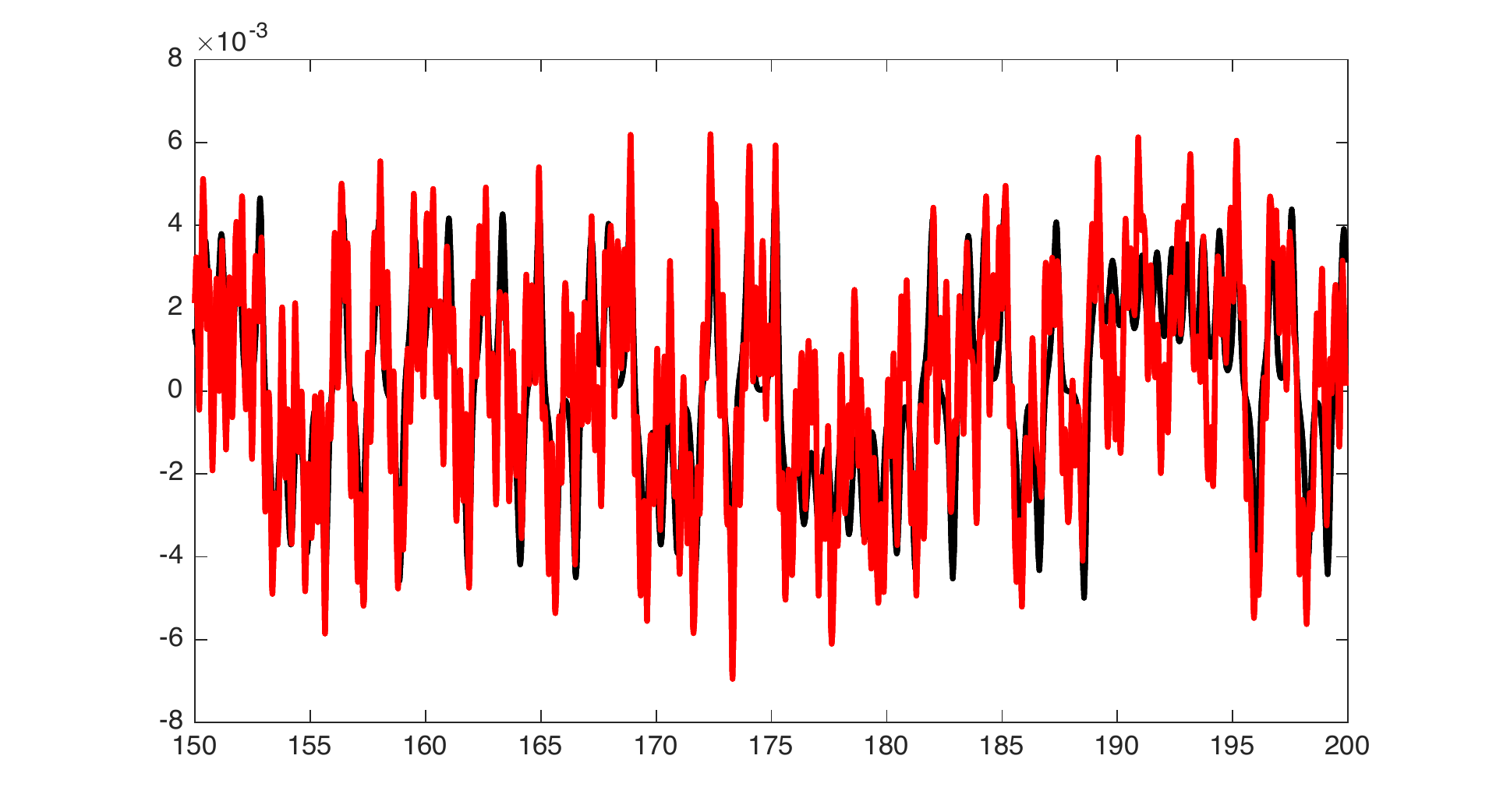}
\put(52,0){$t$}
\put(5.5,30){$v_1$}
\end{overpic}
\caption{Last 50 time unites of response of idealized HAVOK model with integer coefficients, forced by the $v_{15}$ time series from the Lorenz system.}\label{Fig:LorenzP9B}
\end{center}
\end{figure}

This idealized system is forced with the signal $v_{15}$ from the full Lorenz system, and the dynamic response is shown in Fig.~\ref{Fig:LorenzP9}.  
As shown in the zoom-ins in Figs.~\ref{Fig:LorenzP9A} and~\ref{Fig:LorenzP9B}, the agreement is remarkably good for the first 50 time units, although the idealized model performance degrades over time, as shown for the final 50 time units.   

\newpage

The eigenvalues of the full HAVOK model and the integer-valued approximation are provided in Tab.~\ref{Tab:HAVOKEigs}.  
There is reasonably good agreement between the eigenvalues, with the integer-valued model being exactly integrable, so that trajectories reside on a quasi-periodic orbit.  
This is to be expected with a Koopman operator model of a dynamical system.  
Interestingly, a number of the imaginary parts of the eigenvalue pairs are near multiples of each other (e.g., 22.0058 is nearly a multiple of 11.1600).  
This is also to be expected in a Koopman model, as integer multiples of an eigenvalue of the Koopman operator are also eigenvalues corresponding to that integer power of the Koopman eigenfunction.  

\begin{table}[b]
\caption{Eigenvalues of the HAVOK model and integer-valued approximation.}\label{Tab:HAVOKEigs}
\vspace{-.2in}
\begin{center}
\begin{tabular}{|l||l|}
\hline
HAVOK Model & Integer-Valued Model\\\hline\hline
~~~$0.0000 \pm 2.9703i$ & $0.0000 \pm 3.0725i$\\\hline
$-0.0008  \pm 11.0788i$ & $0.0000  \pm 11.1600i$\\\hline
$-0.0035  \pm 21.2670i$ & $0.0000  \pm 22.0058i$\\\hline
$-0.0107  \pm 34.5458i$ & $0.0000  \pm 35.5714i$\\\hline
$-0.0233  \pm 51.4077i$ & $0.0000  \pm 52.3497i$\\\hline
$-0.0301  \pm 72.7789i$ &$0.0000  \pm 73.5112i$ \\\hline
~~~$0.0684  \pm 101.6733i$ &$0.0000  \pm 102.2108i$ \\\hline
\end{tabular}
\end{center}
\end{table}

\newpage
\subsection{Connection to almost-invariant sets and Perron-Frobenius}
The Koopman operator is the dual, or left-adjoint, of the Perron-Frobenius operator, which is also called the transfer operator or the push-forward operator on the space of probability density functions.  
Thus, Koopman analysis typically describes the evolution of measurements from a single trajectory, while Perron-Frobenius analysis describes an ensemble of trajectories.  

Because of the close relationship of the two operators, it is interesting to compare the Koopman (HAVOK) analysis with the almost-invariant sets from the Perron-Frobenius operator.  
Almost-invariant sets represent dynamically isolated phase space regions, in which the trajectory resides for a long time, and with only a weak communication with their surroundings. These sets are almost invariant under the action of the dynamics and are related to dominant eigenvalues and eigenfunctions of the Perron-Frobenius operator. They can be numerically determined from its finite-rank approximation by discretizing the phase space into small boxes and computing a large, but sparse transition probability matrix of how initial conditions in the various boxes flow to other boxes in a fixed amount of time; for this analysis, we use the same $T = 0.1$ used for the length of the $\mathbf{U}$ vectors in the HAVOK analysis. 
Following the approach proposed by \cite{Froyland2005physd},  
almost-invariant sets can then be estimated by computing the associated reversible transition matrix and level-set thresholding its right eigenvectors. 

\begin{figure}[b!]
\begin{center}
\vspace{-.2in}
\begin{overpic}[width=.75\textwidth]{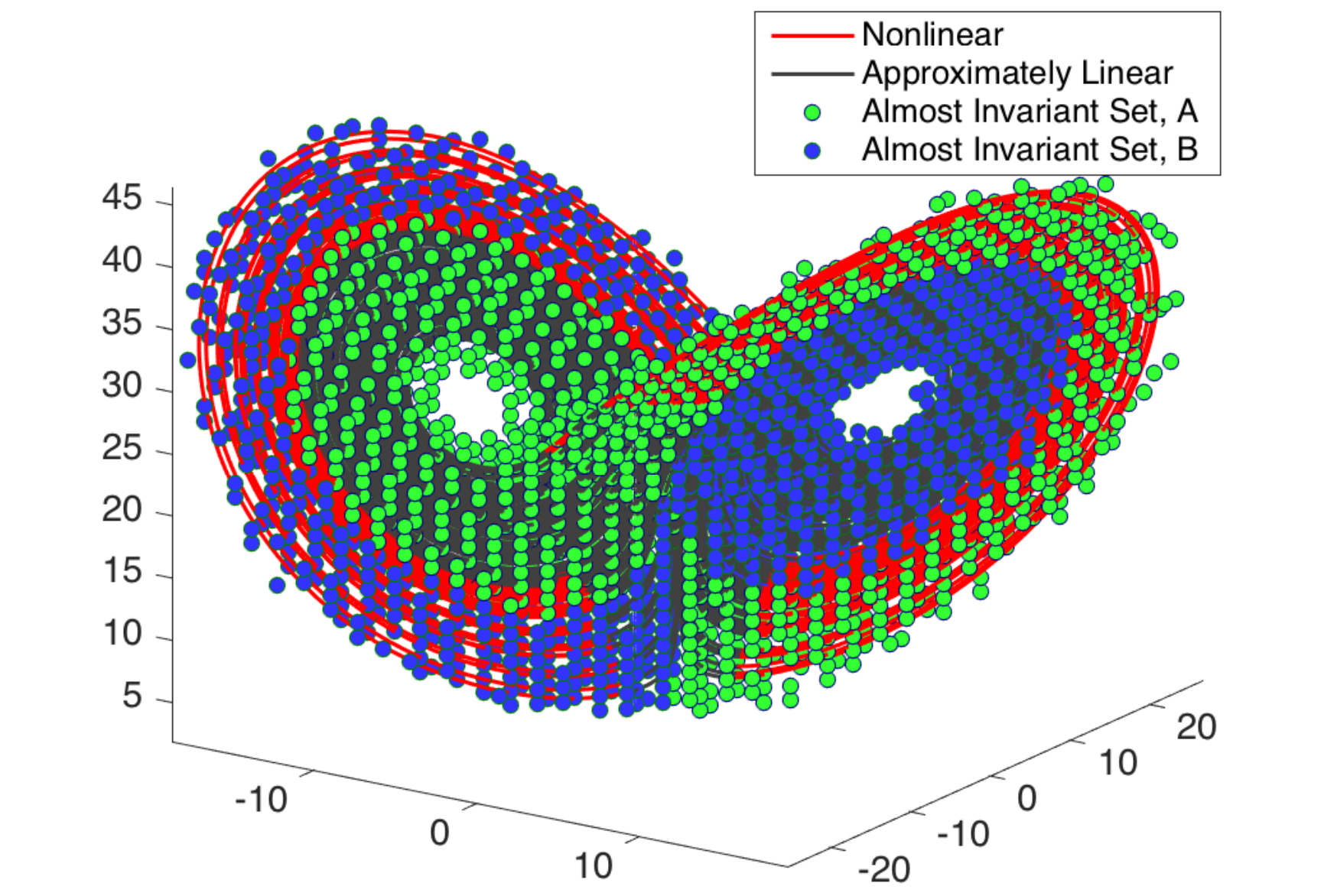}
\put(82,6){$y$}
\put(26,2){$x$}
\put(4,35){$z$}
\end{overpic}
\vspace{-.1in}
\caption{Lorenz attractor visualized using both the HAVOK approximately linear set as well as the Perron-Frobenius almost-invariant sets.}\label{Fig:LorenzPF}
\vspace{-.4in}
\end{center}
\end{figure}

The almost-invariant sets of the Perron-Frobenius operator for the Lorenz system are shown in Fig.~\ref{Fig:LorenzPF}.  
There are two sets, and each set corresponds to the near basin of one of the attractor lobes as well as the outer basin of the opposing attractor lobe and the bundle of trajectories that connect them.  
These two almost-invariant sets dovetail to form the complete Lorenz attractor.  
Underneath the almost-invariant sets, the Lorenz attractor is colored by the thresholded magnitude of the nonlinear forcing term in the HAVOK model, which partitions the attractor into two sets corresponding to regions where the flow is approximately linear (inner black region) and where the flow is strongly nonlinear (outer red region).  
Remarkably, the boundaries of the almost-invariant sets of the Perron-Frobenius operator closely match the boundaries of the linear and nonlinear regions from the HAVOK analysis.  
However, this may not be surprising, as the boundaries of these sets are dynamic separatrices that mark the boundaries of qualitatively different dynamics. 

\newpage
\section{HAVOK analysis on additional chaotic systems}\label{Sec:Results}
\begin{figure*}[b!]
\begin{center}
\begin{overpic}[width=.95\textwidth]{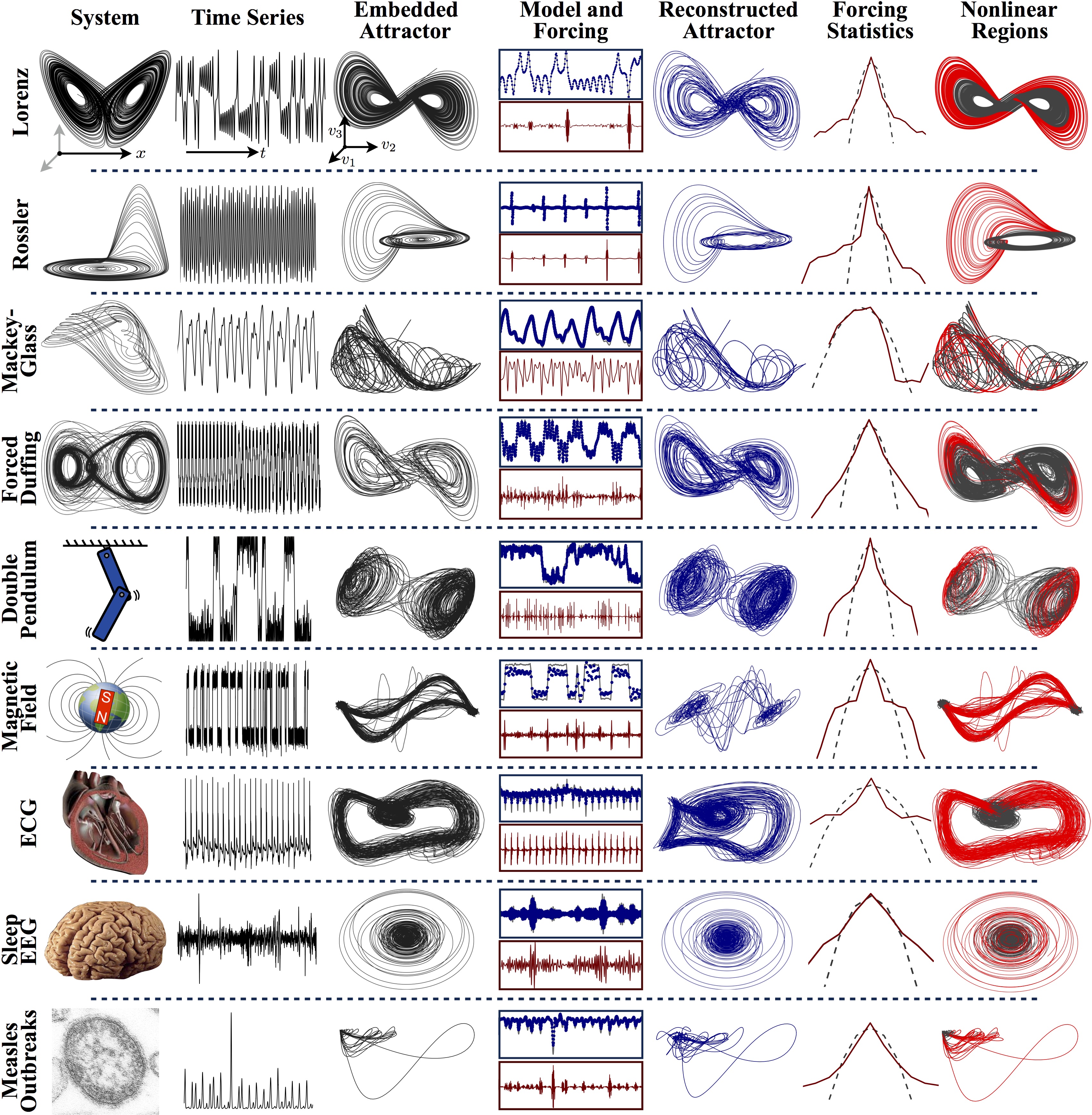}
\end{overpic}
\end{center}
\vspace{-.15in}
\caption{\small Summary of HAVOK analysis applied to numerous examples, including analytical systems (Lorenz and R\"ossler), delay differential equations (Mackey-Glass), stochastic/forced equations (forced Duffing and stochastic magnetic field reversal), and systems characterized from real-world data (electrocardiogram, electroencephalogram, and measles outbreaks).  
The model prediction is extremely accurate for the first five analytical cases, providing faithful attractor reconstruction and predicting dominant transient and intermittent events.  
In the stochastic and data examples, the model predicts attractor topology and switching dynamics.  
In the case of measles outbreaks, the forced linear model predicts large transients corresponding to outbreaks.  
The majority of the examples are characterized by nearly symmetric forcing distributions with long tails (Gaussian forcing shown in black dashed line), corresponding to rare events.  
}\label{Fig:Comprehensive}
\vspace{-.1in}
\end{figure*}

The HAVOK analysis is applied to analytic and real-world systems in Fig.~\ref{Fig:Comprehensive}.  
Code for every example is available at: http://faculty.washington.edu/sbrunton/HAVOK.zip.

The examples span a wide range of systems, including canonical chaotic dynamical systems, such as the Lorenz and R\"{o}ssler systems, 
and the double pendulum, which is among the simplest physical systems that exhibits chaotic motion.  
As a more realistic example, we consider a stochastically driven simulation of the Earth's magnetic field reversal~\cite{Petrelis2009prl}, where complex magnetohydrodynamic (MHD) equations are modeled as a dynamo driven by turbulent fluctuations.   
In this case, the exact form of the attractor is not captured by the linear model, although the attractor switching, corresponding to magnetic field reversal, is preserved.  
In the final three examples, we explore the method on data collected from an electrocardiogram (ECG), electroencephalogram (EEG), and recorded measles cases in New York City over a 36 year timespan from 1928 to 1964.    

In each example, the qualitative attractor dynamics are captured, and importantly, large transients and intermittent phenomena are predicted by the model.  
These large transients and intermittent events correspond to coherent regions in phase space where the forcing is large (right column of Fig.~\ref{Fig:Comprehensive}, red).   
Regions where the forcing is small (black) are well-modeled by a Koopman linear system in delay coordinates.  
Large forcing often precedes intermittent events (lobe switching for Lorenz system and magnetic field reversal, or bursting measles outbreaks), making this signal \emph{predictive}.  

The theory presented here is particularly useful because the forced linear models generalize beyond the training data on which they are constructed.  
It is possible to test the HAVOK model on data from a test trajectory that was not used to train the model.  
The forcing signal $v_r(t)$ is extracted from the new validation time series and fed as an input into the model in \eqref{Eq:ChaosModel}, resulting in a faithful reproduction of the attractor dynamics, including lobe switching events. 
The forcing term $v_r(t)$ is obtained from a relatively short window of time compared with a switching event, so it may be measured and used for prediction in realistic applications.  

\subsection{Parameters}
All systems in this paper are described in this section, including details about the equations of motion, data collection, and key characteristics of each system.  The simulation parameters for the numerical examples are in Tab.~\ref{Table:Parameters}, and the parameters used for the HAVOK analysis are in Tab.~\ref{Table:HAVOK}.  

The example systems range from ordinary differential equation (ODE) models to delay differential equations (DDE) and stochastic differential equations (SDE), as well as systems characterized purely by measurement data.  
We also consider the forced Duffing oscillator, which is not chaotic for the parameters examined, but provides a weakly nonlinear example to investigate.

\begin{table}[where!]
\caption{Parameters used to simulate the various numerical example systems.}\label{Table:Parameters}
\vspace{-.2in}
\begin{center}
\begin{tabular}{|c||c|c|}
\hline
\bf System & \bf Parameters & \bf Initial Conditions \\\hline\hline
\bf Lorenz  & $\sigma=10, \rho=28, \beta=8/3$& $(-8, 8, 27)$ \\\hline
\bf R\"{o}ssler  & $a=0.1,b=0.1,c=14$ & $(1,1,1)$ \\\hline
\bf Mackey-Glass  & $\beta=2,\tau=2,n=9.65,\gamma =1$ & $0.5$ \\\hline
\bf Unforced Duffing & $\delta=0, \alpha=-1, \beta=5, \gamma=0, \omega=0$ & Multiple  \\\hline
\bf Forced Duffing  & $\delta=0.02, \alpha=1, \beta=5, \gamma=8, \omega=0.5$ & Multiple\\\hline
\bf Double Pendulum  & $l_1=l_2=m_1=m_2=1$, $g=10$ & $(\pi/2,\pi/2,-0.01,-0.005)$ \\\hline
\bf Magnetic Field  & $B_1=-0.4605i, B_2=-1+0.12i$, & $0.1$\\
&  $B_3=0.4395i$, $B_4=-0.06 -0.12i$, & \\
 & $b_1=b_4=0.25$, and $b_2=b_3=0.07$,  $\mu= \nu=0$ & \\\hline
\end{tabular}
\end{center}
\vspace{-.25in}
\end{table}

\begin{table}[where!]
\caption{Summary of systems and HAVOK analysis parameters for each example.  *For the measles data, the original 431 data samples are spaced 2 weeks apart; this data is interpolated at 0.2 week resolution, and the new data is stacked into $q=50$ rows.}\label{Table:HAVOK}
\vspace{-.2in}
\small
\begin{center}
\begin{tabular}{|c||c|c|c|c|c|c|c|}
\hline
\bf System &\bf Type & \bf Measured & \bf\# Samples & \bf Time step&\bf \# Rows in $\mathbf{H}$&\bf Rank &  \bf Energy in $r$ \\
&&\bf Variable &  $m$ & $\Delta t$ & $q$ & $r$ &\bf  modes  (\%)  \\\hline\hline
\bf Lorenz  & ODE& $x(t)$& 200,000& 0.001&100 & 15  & 100 \\\hline
\bf R\"{o}ssler  & ODE&$x(t)$& 500,000 & 0.001& 100& 6 & 99.9999997\\\hline
\bf Mackey-Glass  & DDE&$x(t)$& 100,000& 0.001 &100 &4 &99.9999 \\\hline
\bf Unforced Duffing &ODE &$x(t)$& 360,000&0.001 &100 &5 & 97.1\\\hline
\bf Forced Duffing  & ODE&$x(t)$& 1,000,000&0.001 & 100 &5 &99.9998 \\\hline
\bf Double Pendulum  & ODE&$\sin(\theta_1(t))$& 250,000&0.001 & 100 & 5&99.997 \\\hline
\bf Magnetic Field  & SDE& $\text{Re}(A)$& 100,000& 1 year & 100& 4& 55.2\\\hline
\bf ECG  & Data& Voltage& 45,000& 0.004 s& 25& 5&67.4 \\\hline
\bf Sleep EEG  & Data& Voltage&100,000 & 0.01 s&1,000 & 4 & 7.7\\\hline
\bf Measles Outbreak  &Data & Cases & 431& 2 weeks& $~~$50$^*$& 9& 99.78\\\hline
\end{tabular}
\end{center}
\vspace{-.25in}
\end{table}

\subsection{Description of all examples}
\subsubsection{Duffing oscillator}\label{Sec:Res:Duffing}
The Duffing oscillator is a simple dynamical system with two potential wells separated by a saddle point, as shown in the middle panel of Fig.~\ref{fig03}.
The dynamics are given by
\begin{eqnarray}
\label{Eq:Duffing}
\ddot{x} +\delta \dot{x} + \alpha x + \beta x^3 = \gamma\cos(\omega t).
\end{eqnarray}
It is possible to suspend variables to obtain a coupled system of first-order equations:
\begin{subequations}
\label{Eq:DuffingSys}
\begin{eqnarray}
\dot{x} & = & v\\
\dot{v} & = & -\delta v - \alpha x - \beta x^3 + \gamma \cos(\omega t).
\end{eqnarray}
\end{subequations}

\subsubsection{Chaotic Lorenz system}\label{Sec:Res:Lorenz}
The Lorenz system~\cite{Lorenz1963jas} is a canonical model for chaotic dynamics:  
\begin{subequations}
\label{Eq:Lorenz2}
\begin{eqnarray}
\dot{x} & = & \sigma (y - x)\\
\dot{y} & = & x(\rho -z) - y\\
\dot{z} & = & x y - \beta z.
\end{eqnarray}
\end{subequations}

\subsubsection{R\"ossler system}\label{Sec:Res:Rossler}
The R\"ossler system is given by
\begin{subequations}
\label{Eq:Rossler}
\begin{eqnarray}
\dot{x} & = & -y-z\\
\dot{y} & = & x+ay\\
\dot{z} & = & b + z(x-c)
\label{Eq:Rossler}
\end{eqnarray}
\end{subequations}
This system exhibits interesting dynamics, whereby the trajectory is characterized by oscillation in the $x-y$ plane punctuated by occasional transient growth and decay in the $z$ direction.  
We refer to this behavior as \emph{bursting}, as it corresponds to a transient from an attractor to itself, as opposed to the lobe switching between attractor lobes in the Lorenz system.  

\subsubsection{Mackey-Glass delay differential equation}\label{Sec:Res:MG}
The Mackey-Glass equation is a canonical example of a delay differential equation, given by
\begin{eqnarray}
\label{Eq:MackeyGlass}
\dot{x}(t) = \beta \frac{x(t-\tau)}{1+x(t-\tau)^n} - \gamma x(t),
\end{eqnarray}
with $\beta=2$, $\tau=2$, $n=9.65$, and $\gamma=1$.  
The current time dynamics depend on the state $x(t-\tau)$ at a previous time $\tau$ in the past.  

\subsubsection{Double pendulum}\label{Sec:Res:DP}
The double pendulum is among the simplest physical systems that exhibits chaos.  The numerical simulation of the double pendulum is sensitive, and we implement a variational integrator based on the Euler-Lagrange equations
\begin{eqnarray}
\frac{d}{dt}\frac{\partial L}{\partial \dot{\mathbf{q}}} - \frac{\partial L}{\partial \mathbf{q}}=0
\end{eqnarray}
where the Lagrangian $L=T-V$ is the kinetic ($T$) minus potential ($V$) energy.  For the double pendulum, $\mathbf{q}=\begin{bmatrix}\theta_1 & \theta_2\end{bmatrix}^T$, and the Lagrangian becomes:
\begin{align*}
L = T-V =& \frac{1}{2}(m_1+m_2)l_1\dot{\theta}_1^2 + \frac{1}{2}m_2l_2^2\dot{\theta}_2^2 + m_2l_1l_2\dot{\theta}_1\dot{\theta}_2\cos(\theta_1-\theta_2)\\
& -(m_1+m_2)l_1g(1-\cos(\theta_1))-m_2l_2g\left(1-\cos(\theta_2)\right).
\end{align*}
We integrate the equations of motion with a variational integrator derived using a trapezoidal approximation to the action integral:
\begin{align*}
\delta \int_a^b L(\mathbf{q},\dot{\mathbf{q}},t)dt = 0.
\end{align*}

Because the mean of $\theta_1$ drifts after a revolution, we use $x(t)=\cos(2\theta_1(t))$ as a measurement.  

\subsubsection{Earth's magnetic field reversal}\label{Sec:Res:Magnetic}
The Earth's magnetic field is known to reverse over geological time scales~\cite{Cox1964science,Guyodo1999nature}.  
Understanding and predicting these rare events is an important challenge in modern geophysics.  
It is possible to model the underlying dynamics that give rise to magnetic field switching by considering the turbulent magnetohydrodynamics inside the Earth.  
A simplified model~\cite{Petrelis2008jp,Petrelis2009prl,Petrelis2010ptrsla} may be obtained by modeling the turbulent fluctuations as stochastic forcing on a dynamo.

This is modeled by the following differential equation in terms of the magnetic field $A$:
\begin{align}
\frac{d}{dt}A = \mu A + \nu \bar{A} + B_1A^3 + B_2A^2\bar{A} + B_3A\bar{A}^2+B_4\bar{A}^3 + f
\end{align}
where $\bar{A}$ is the complex conjugate of $A$ and the stochastic forcing $f$ is given by
\begin{align*}
f = \left(b_1\xi_1 + ib_3\xi_3\right)\text{Re}(A) + \left(b_2\xi_2+ib_4\xi_4\right)\text{Im}(A).
\end{align*}
The variables $\xi$ are Gaussian random variables with a standard deviation of $5$.  The other parameters are given by $\mu=1$, $\nu=0$, $B_1=-0.4605i, B_2=-1+0.12i, B_3=0.4395i$, $B_4=-0.06 -0.12i$, $b_1=b_4=0.25$, and $b2=b_3=0.07$.

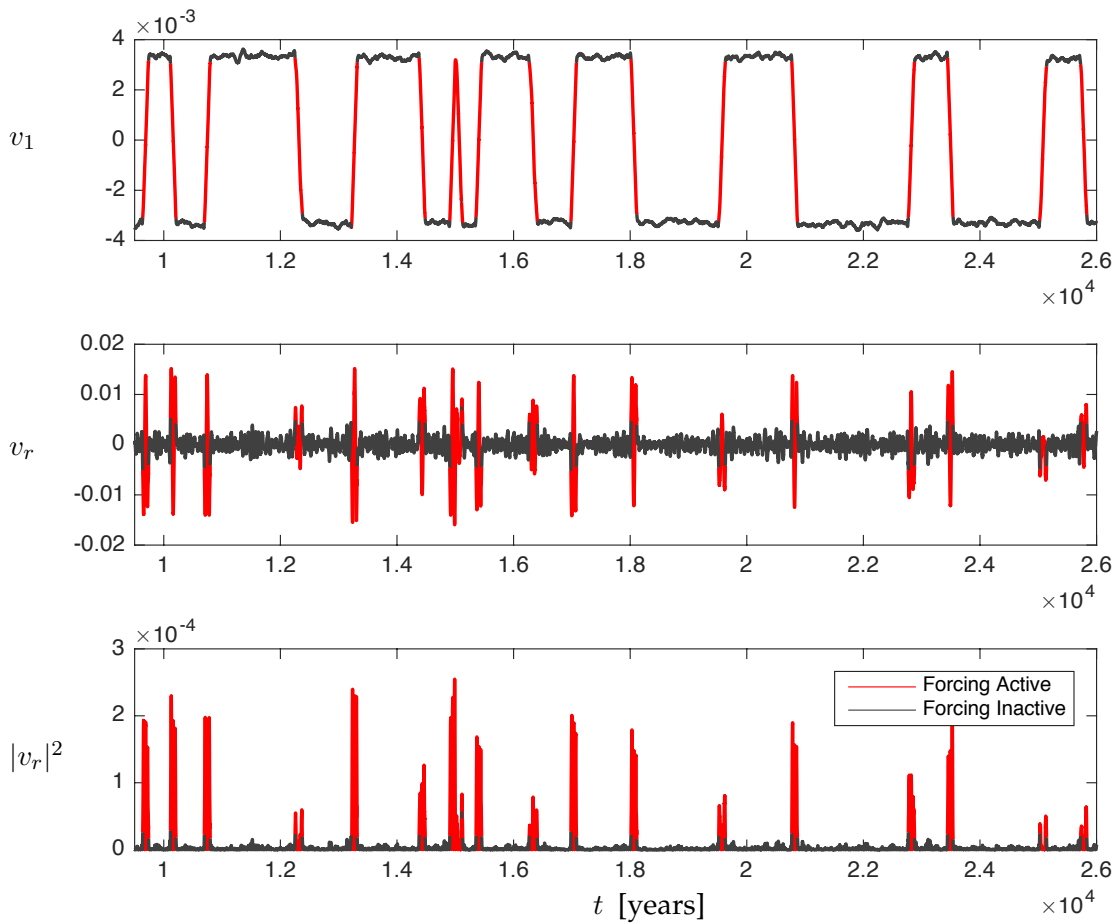
\begin{figure}
\begin{center}
\begin{overpic}[width=\textwidth]{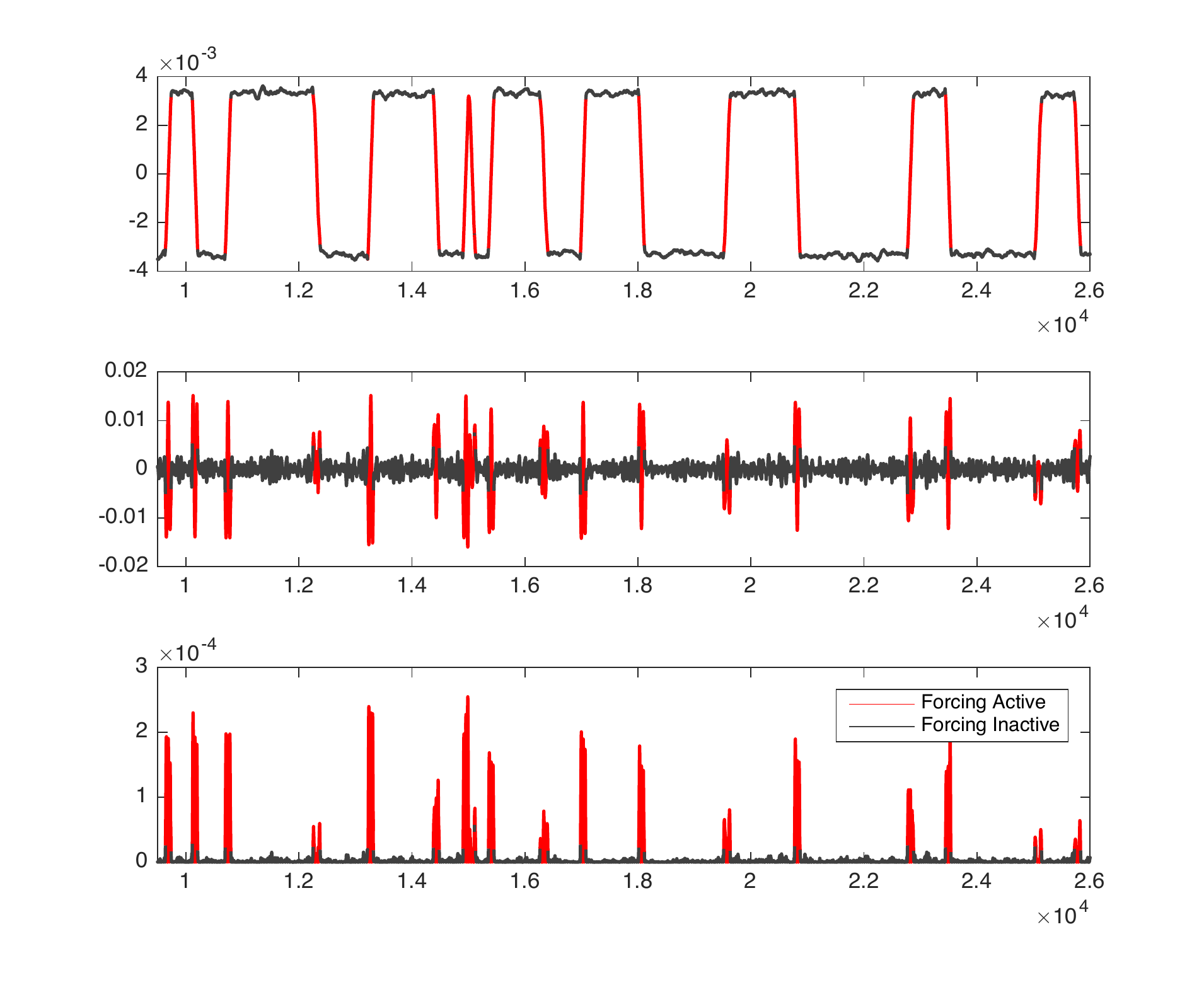}
\put(50,5){$t~$ [years]}
\put(3,17){$|v_{r}|^2$}
\put(3,42){$v_{r}$}
\put(3,67){$v_1$}
\end{overpic}
\vspace{-.45in}
\caption{Prediction of Earth's magnetic field reversal using the forcing term $v_r$.}\label{Fig:MHD}
\vspace{-.25in}
\end{center}
\end{figure}

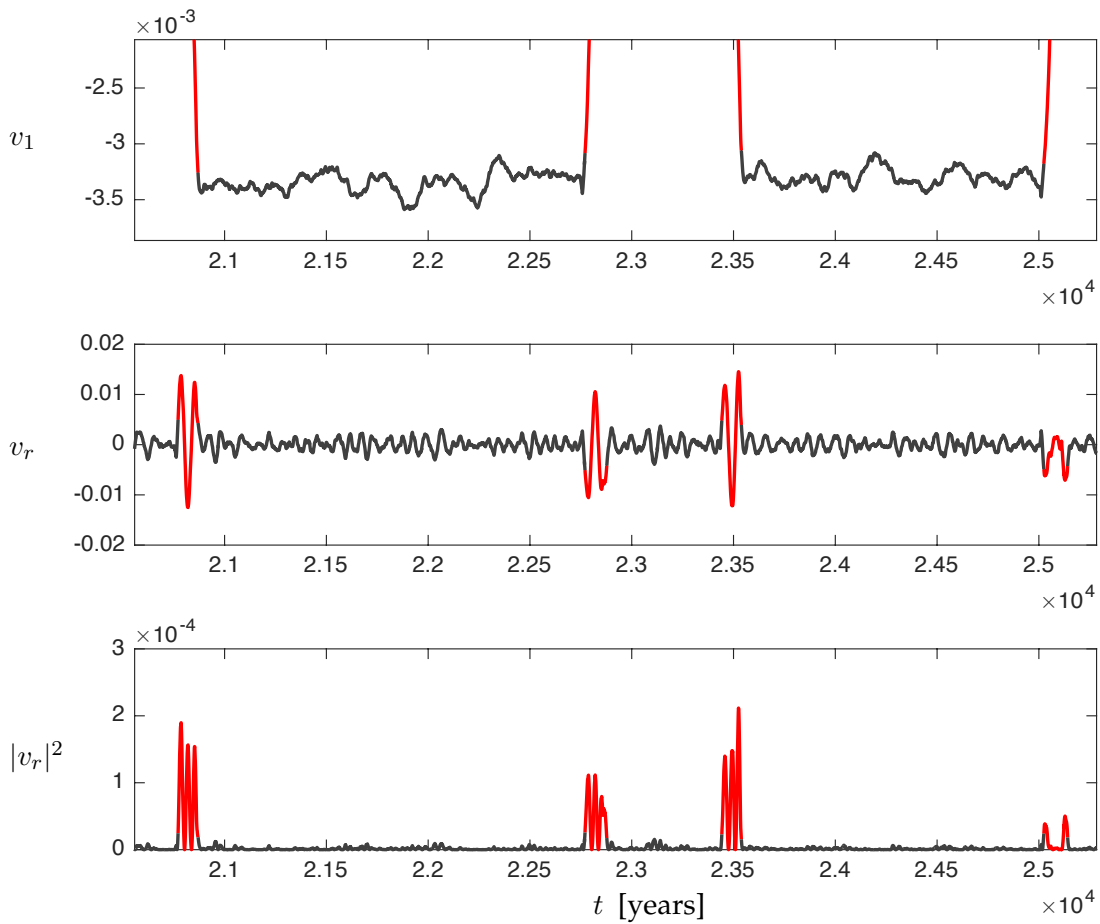
\begin{figure}
\begin{center}
\begin{overpic}[width=\textwidth]{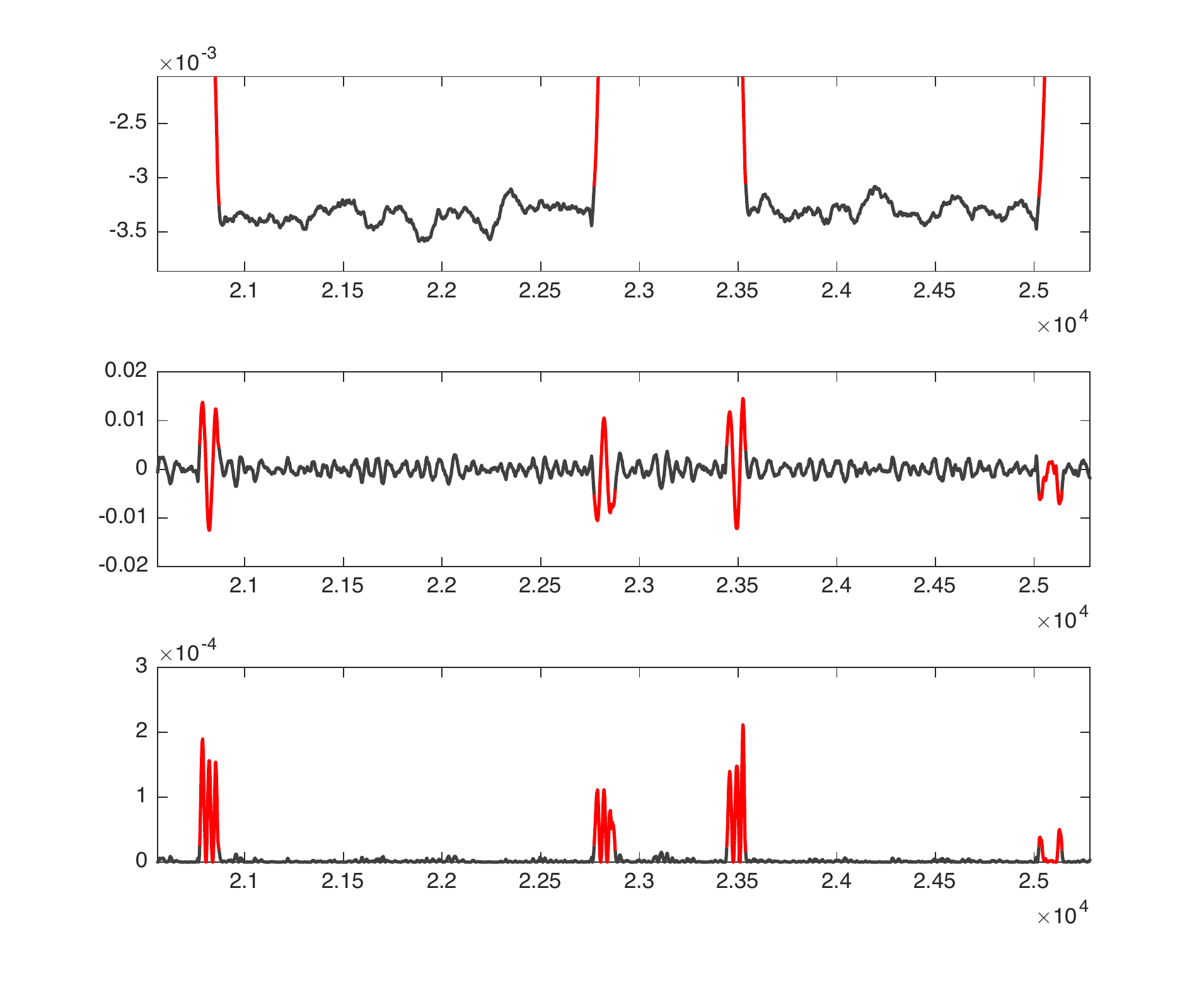}
\put(50,5){$t~$ [years]}
\put(3,17){$|v_{r}|^2$}
\put(3,42){$v_{r}$}
\put(3,67){$v_1$}
\end{overpic}
\vspace{-.45in}
\caption{Prediction of Earth's magnetic field reversal using the forcing term $v_r$. Zoom-in of field reversal events in Fig.~\ref{Fig:MHDZOOM}.}\label{Fig:MHDZOOM}
\end{center}
\end{figure}

Figure~\ref{Fig:MHD} shows the $v_1$ and $v_r$ eigen time series for the Earth magnetic field reversal data.  
The signal in $v_1$ accurately represents the global field switching dynamics.  
Color-coding the trajectory by the magnitude of the forcing in $v_r$, it is clear that the field reversal events correspond to large forcing.  
Analyzing two of these field switching events more closely in Fig.~\ref{Fig:MHDZOOM}, it appears that the forcing signal becomes active before the signal in $v_1$ exceeds the expected variance.  
This means that the forcing signal provides a prediction of field reversal before there is a clear statistical signature in the dominant $v_1$ time series.

\subsubsection{Electrocardiogram (ECG)}\label{Sec:Res:ECG}
An electrocardiogram (ECG) measures the electrical activity of the heart, producing the characteristic spiking pulses associated with each heartbeat.  
Data from ECG has long been analyzed using delay embeddings to quantify the fractal dimension of ECG recordings, among other quantities~\cite{Schreiber1996chaos,Richter1998pre}.  
In this analysis, we use the ECG signal \texttt{qtdb/sel102} that was used in~\cite{Keogh2005icdm}, adapted from the PhysioNet database~\cite{Laguna1997cc,Goldberger2000circulation}.  
This signal corresponds to $T=380$ seconds of data with a sampling rate of $250$ Hz.  

\subsubsection{Electroencephalogram (EEG)}\label{Sec:Res:EEG}
An electroencephalogram (EEG) is a nonintrusive measurement of brain activity achieved through electrodes placed on the scalp.  
A time series of voltage is recorded from these electrodes, and although typically the spectral content of EEG is analyzed, there are also numerous examples of time series EEG analysis~\cite{Pritchard1992measuring,Acharya2005non,Kannathal2005entropies,Stam2005nonlinear,lainscsek2015delay}.  
Although it is possible to obtain vast quantities of data, possibly using an array of electrodes, the voltage signal is only a rough proxy for brain activity, as signals must pass through thick layers of dura, cerebrospinal fluid, skull, and scalp.  

Data is available at~\cite{Goldberger2000circulation,Kemp2000ieeetbe}: https://physionet.org/pn4/sleep-edfx/

\subsubsection{Measles outbreaks}\label{Sec:Res:Measles}
Time series analysis has long been applied to understand and model disease epidemics.  
Measles outbreaks have been particularly well studied, partially because of the wealth of accurate historical data.  
The seminal work of Sugihara and May~\cite{Sugihara1990nature} use Takens embedding theory to analyze Measles outbreaks.  
It was also shown that Measles outbreaks are chaotic~\cite{Schaffer1985bs}. 
Here, we use data for Measles outbreaks in New York City (NYC) from 1928 to 1964, binned every 2 weeks~\cite{London1973aje}.

\begin{figure}[where!]
\begin{center}
\begin{overpic}[width=.95\textwidth]{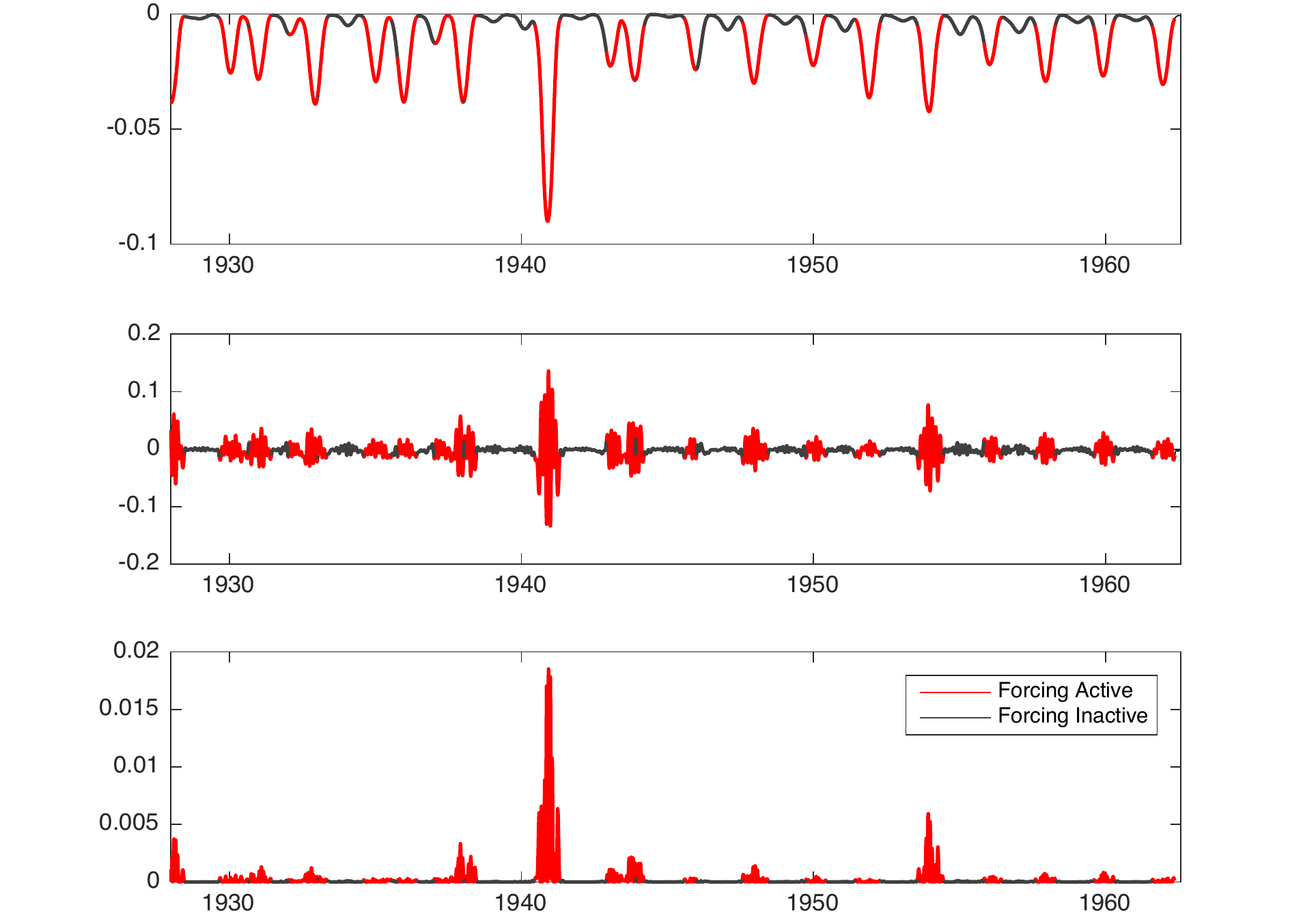}
\put(55,-1){$t$}
\put(3,14){$|v_{r}|^2$}
\put(3,36){$v_{r}$}
\put(3,60){$v_1$}
\end{overpic}
\vspace{-.05in}
\caption{Prediction of Measles outbreaks using the forcing term $v_r$.}\label{Fig:Measles}
\vspace{-.3in}
\end{center}
\end{figure}
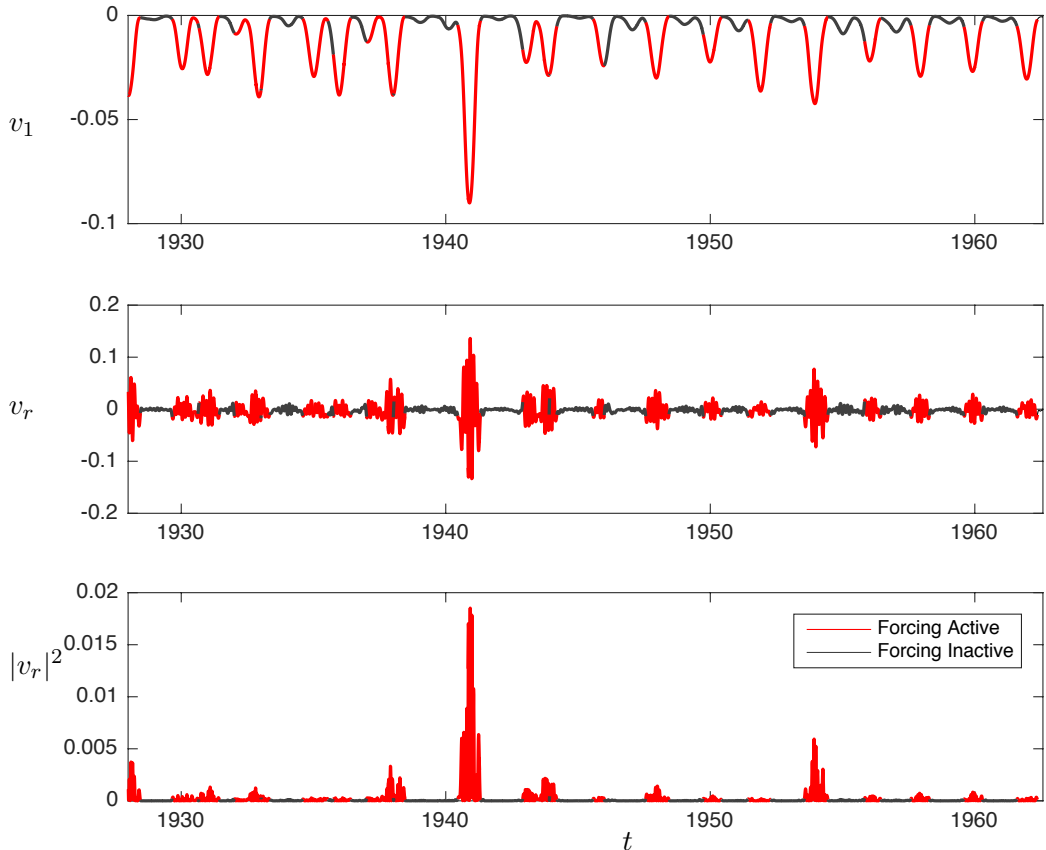

Figure~\ref{Fig:Measles} shows the $v_1$ and $v_r$ eigen time series for the Measles outbreaks data.  
The $v_1$ signal provides a signature for the severity of the outbreak, with larger negative values corresponding to more cases of Measles.  
The forcing signal accurately captures many of the outbreaks, most notably the largest outbreak after 1940.  
This outbreak was preceded by two small dips in $v_1$, which may have resulted in false positives using other prediction methods.  
However, the forcing signal becomes large directly preceding the outbreak. 

\section{Discussion}
In summary, we have presented a data-driven procedure, the HAVOK analysis, to identify an intermittently forced linear system representation of chaos.  
This procedure is based on machine learning regressions, Takens' embedding, and Koopman operator theory.  
The activity of the forcing signal in these models predicts intermittent transient events, such as lobe switching and bursting, and partitions phase space into coherent linear and nonlinear regions.  
More generally, the search for intrinsic or natural measurement coordinates is of central importance in finding simple representations of complex systems, and this will only become increasingly important with growing data.  
Simple, linear representations of complex systems is a long sought goal, providing the hope for a general theory of nonlinear estimation, prediction, and control.  
This analysis will hopefully motivate novel strategies to measure, understand, and control~\cite{Shinbrot1993nature} chaotic systems in a variety of scientific and engineering applications.  

\newpage
\bibliographystyle{plain}
\bibliography{references}

\newpage
\begin{table}
\begin{center}
\vspace{-.15in}
\includegraphics[width=.65\textwidth]{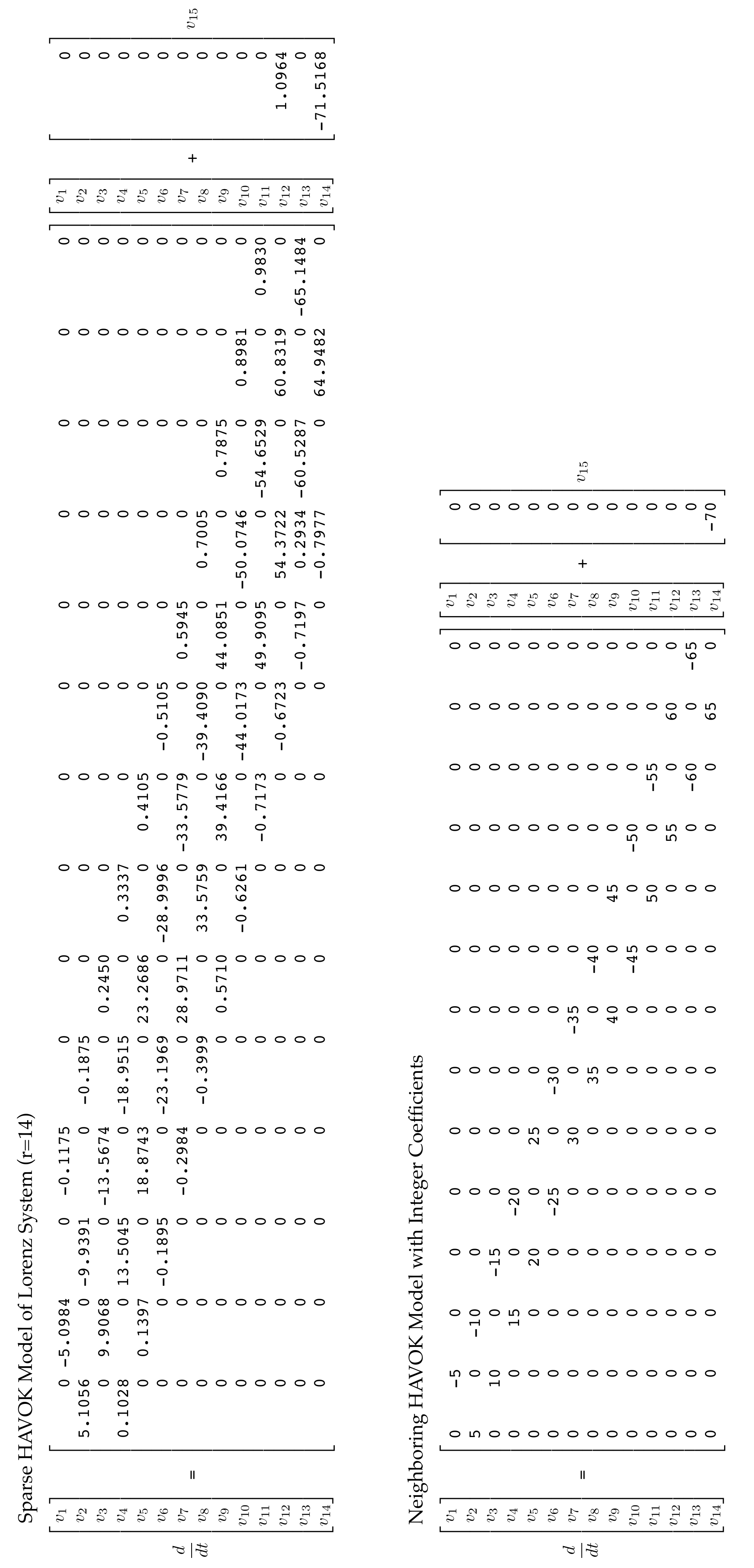}
\caption{HAVOK model for Lorenz system.}\label{Fig:HAVOKLorenzExact}
\end{center}
\end{table}

\end{document}